\documentclass{article}
\usepackage{tikz,tikz-qtree,tikz-qtree-compat,adjustbox}
\usepackage{etoolbox}
\usepackage{url}
\usepackage{amssymb}
\usepackage{amsmath}
\usepackage{graphics}
\usepackage[utf8]{inputenc}
\usepackage[mathscr]{eucal}
\usepackage{amsthm,amssymb}
\newtheorem{Theorem}{Theorem}

\newtheorem{Lemma}[Theorem]{Lemma}
\newtheorem{Corollary}[Theorem]{Corollary}

\newtheorem{Conjecture}[Theorem]{Conjecture}
\newtheorem{Remark}[Theorem]{Remark}
\newtheorem*{Definition}{Definition}
\newenvironment{keywords}{{\bf Keywords: }}{}
\usepackage{hyperref}
\usepackage{subfigure}
\usepackage{longtable}

\usepackage{tikz,tikz-qtree,tikz-qtree-compat,adjustbox}

\def\maxomega{15}

\renewcommand\leq\leqslant
\renewcommand\geq\geqslant

\newcommand{\circledgenerator}[2][]{\tikz[baseline=(char.base)]{\node[shape = circle, draw, inner sep = 2pt,fill=black!40](char) {\phantom{\ifblank{#1}{#2}{#1}}};\node at (char.center) {\makebox(0,0){#2}};}}
\newcommand{\circlednongap}[2][]{\tikz[baseline=(char.base)]{\node[shape = circle, draw, inner sep = 2pt,fill=black!20](char) {\phantom{\ifblank{#1}{#2}{#1}}};\node at (char.center) {\makebox(0,0){#2}};}}
\newcommand{\circledgap}[2][]{\tikz[baseline=(char.base)]{\node[shape = circle, draw, inner sep = 2pt](char) {\phantom{\ifblank{#1}{#2}{#1}}};\node at (char.center) {\makebox(0,0){#2}};}}
\robustify{\circlednongap}
\robustify{\circledgap}
\robustify{\circledgenerator}
\newcommand{\gap}[1]{\scalebox{.5}{\circledgap[00]{\large 
      {\color{black!30} \sf #1}}}}
\newcommand{\nongap}[1]{\scalebox{.5}{\circlednongap[00]{\large\sf #1}}}
\newcommand{\generator}[1]{\scalebox{.5}{\circledgenerator[00]{\large\sf #1}}}

\title{Quasi-ordinarization transform \\of a numerical semigroup}
\author{Maria Bras-Amorós, Hebert Pérez-Rosés, José Miguel Serradilla Merinero
\thanks{This work was partly supported by the Catalan Government under grant 2017 SGR 00705 and by the Spanish Ministry of Economy and Competitivity under grant TIN2016-80250-R. }}

\begin{document}

\maketitle

\begin{abstract}In this study, we present the notion of the quasi-ordinarization transform of a numerical semigroup. The set of all semigroups of a fixed genus can be organized in a forest whose roots are all the quasi-ordinary semigroups of the same genus. This way, we approach the conjecture on the increasingness
of the cardinalities of the sets of numerical semigroups of each given genus.
We analyze the number of nodes at each depth in the forest and propose new conjectures.
Some properties of the quasi-ordinarization transform are presented, as well as some relations between the ordinarization and quasi-ordinarization transforms.
\end{abstract}


\section{Introduction}
A numerical semigroup is a cofinite submonoid of ${\mathbb N}_0$ under addition,
where ${\mathbb N}_0$ is the set of nonnegative integers.

While the symmetry of structures has traditionally been studied with the aid of groups, it is also possible to relax the definition of symmetry, so as to describe some forms of symmetry that arise in quasicrystals, fractals, and other natural phenomena, with the aid of semigroups or monoids, rather than groups. For example, Rosenfeld and Nordahl \cite{Ros16} lay the groundwork for such a theory of symmetry based on semigroups and monoids, and they cite some applications in Chemistry. 

Suppose that $\Lambda$ is a numerical semigroup. The elements in the complement ${\mathbb N}_0\setminus\Lambda$ are called the {\em gaps} of the semigroup and the number of gaps is its {\em genus}. The {\em Frobenius number} is the largest gap and the {\em conductor} is the non-gap that equals the Frobenius number plus one.
The first non-zero non-gap of a numerical semigroup (usually denoted by $m$) is called its {\em multiplicity}.
An {\em ordinary} semigroup is a numerical semigroup different than ${\mathbb N}_0$ in which all gaps are in a row.
The non-zero non-gaps of a numerical semigroup that are not the result of the sum of two smaller non-gaps
are called the {\em generators} of the numerical semigroup. It is easy to deduce that the set of generators of a numerical semigroup must be co-prime. One general reference for numerical semigroups is \cite{NS}.

To illustrate all these definitions, consider the well-tempered harmonic semigroup
$H=\{0,12,19,24,28,31,34,36,38,40,42,43,45,46,47,48,\dots\}$,
where we use $"\dots"$ to indicate that the semigroup contains consecutively all the integers from the number that precedes the ellipsis.
The semigroup $H$ arises in the mathematical theory of music \cite{Bras:mathematicsandmusic}. 
It is obviously cofinite and it contains zero. One can also check that it is closed under addition. Hence,
it is a numerical semigroup. Its Frobenius number is $44$, its conductor is $45$, its genus is $33$, and its multiplicity is $12$. Its generators are $\{12,19,28,34,42,45,49,51\}$.

The number of numerical semigroups of genus $g$ is denoted $n_g$.
It was conjectured in \cite{Bras:Fibonacci} that the sequence $n_g$ asymptotically behaves as the Fibonacci numbers.
In particular, it was conjectured that each term in the sequence is larger than the sum of the two previous terms, that is, $n_g\geq n_{g-1}+n_{g-2}$ for $g\geq 2$, being each term more and more similar to the sum of the two previous terms as $g$ approaches infinity, more precisely $\lim_{g\to\infty}\frac{n_g}{n_{g-1}+n_{g-2}}=1$ and, equivalently, $\lim_{g\to\infty}\frac{n_g}{n_{g-1}}=\phi=\frac{1+\sqrt{5}}{2}$.
A number of papers deal with the sequence $n_g$ 
\cite{Komeda89,Komeda98,bounds,Anna,Stas,sergi,Zhao,BlancoGarciaPuerto,Kaplan,blancorosales,Bras:ordinarization,odorney,bernardini,fromentin,rgd,kaplanmonthly}.
Alex Zhai proved the asymptotic Fibonacci-like behavior of $n_g$ \cite{Zhai}.
However, it remains not proved that $n_g$ is increasing. This was already conjectured by Bras-Amorós in \cite{ngu}. More information on $n_g$, as well as the list of the first $73$ terms can be found in entry \href{http://oeis.org/A007323}{A007323} of The On-Line Encyclopedia of Integer Sequences \cite{oeis}.

\begin{figure}
  \resizebox{\textwidth}{!}{\begin{tikzpicture}[grow=right, every node/.style = {align=left}]\tikzset{level 1/.style={level distance=3cm}}\tikzset{level 2/.style={level distance=5cm}}\tikzset{level 3/.style={level distance=6cm}}\tikzset{level 4/.style={level distance=8cm}}\tikzset{level 5+/.style={level distance=10cm}}\Tree[.{$\nongap{0}\generator{1}\nongap{2} \dots$} [.{$\nongap{0}\gap{1}\generator{2} \generator{3} \dots$} [.{$\nongap{0}\gap{1}\gap{2}\generator{3} \generator{4} \generator{5} \dots$} [.{$\nongap{0}\gap{1}\gap{2}\gap{3}\generator{4} \generator{5} \generator{6} \generator{7} \dots$} [.{$\nongap{0}\gap{1}\gap{2}\gap{3}\gap{4}\generator{5} \generator{6} \generator{7} \generator{8} \generator{9} \dots$} [.{$\nongap{0}\gap{1}\gap{2}\gap{3}\gap{4}\gap{5}\generator{6} \generator{7} \generator{8} \generator{9} \generator{10} \generator{11} \dots$} [.{$\nongap{0}\gap{1}\gap{2}\gap{3}\gap{4}\gap{5}\gap{6}\generator{7} \generator{8} \generator{9} \generator{10} \generator{11} \generator{12} \generator{13} \dots$}  ]
[.{$\nongap{0}\gap{1}\gap{2}\gap{3}\gap{4}\gap{5}\nongap{6}\gap{7}\generator{8} \generator{9} \generator{10} \generator{11} \nongap{12}\generator{13} \dots$}  ]
[.{$\nongap{0}\gap{1}\gap{2}\gap{3}\gap{4}\gap{5}\nongap{6}\nongap{7}\gap{8}\generator{9} \generator{10} \generator{11} \nongap{12}\nongap{13}\dots$}  ]
[.{$\nongap{0}\gap{1}\gap{2}\gap{3}\gap{4}\gap{5}\nongap{6}\nongap{7}\nongap{8}\gap{9}\generator{10} \generator{11} \nongap{12}\nongap{13}\dots$}  ]
[.{$\nongap{0}\gap{1}\gap{2}\gap{3}\gap{4}\gap{5}\nongap{6}\nongap{7}\nongap{8}\nongap{9}\gap{10}\generator{11} \nongap{12}\nongap{13}\dots$}  ]
[.{$\nongap{0}\gap{1}\gap{2}\gap{3}\gap{4}\gap{5}\nongap{6}\nongap{7}\nongap{8}\nongap{9}\nongap{10}\gap{11}\nongap{12} \nongap{13}\dots$}  ]
 ]
[.{$\nongap{0}\gap{1}\gap{2}\gap{3}\gap{4}\nongap{5}\gap{6}\generator{7} \generator{8} \generator{9} \nongap{10}\generator{11} \dots$} [.{$\nongap{0}\gap{1}\gap{2}\gap{3}\gap{4}\nongap{5}\gap{6}\gap{7}\generator{8} \generator{9} \nongap{10}\generator{11} \generator{12} \nongap{13}\dots$}  ]
[.{$\nongap{0}\gap{1}\gap{2}\gap{3}\gap{4}\nongap{5}\gap{6}\nongap{7}\gap{8}\generator{9} \nongap{10}\generator{11} \nongap{12}\generator{13} \dots$}  ]
[.{$\nongap{0}\gap{1}\gap{2}\gap{3}\gap{4}\nongap{5}\gap{6}\nongap{7}\nongap{8}\gap{9}\nongap{10} \generator{11} \nongap{12}\nongap{13}\dots$}  ]
[.{$\nongap{0}\gap{1}\gap{2}\gap{3}\gap{4}\nongap{5}\gap{6}\nongap{7}\nongap{8}\nongap{9}\nongap{10}\gap{11}\nongap{12} \nongap{13}\dots$}  ]
 ]
[.{$\nongap{0}\gap{1}\gap{2}\gap{3}\gap{4}\nongap{5}\nongap{6}\gap{7}\generator{8} \generator{9} \nongap{10}\nongap{11}\dots$} [.{$\nongap{0}\gap{1}\gap{2}\gap{3}\gap{4}\nongap{5}\nongap{6}\gap{7}\gap{8}\generator{9} \nongap{10}\nongap{11}\nongap{12}\generator{13} \dots$}  ]
[.{$\nongap{0}\gap{1}\gap{2}\gap{3}\gap{4}\nongap{5}\nongap{6}\gap{7}\nongap{8}\gap{9}\nongap{10} \nongap{11}\nongap{12}\nongap{13}\dots$}  ]
 ]
[.{$\nongap{0}\gap{1}\gap{2}\gap{3}\gap{4}\nongap{5}\nongap{6}\nongap{7}\gap{8}\generator{9} \nongap{10}\nongap{11}\dots$} [.{$\nongap{0}\gap{1}\gap{2}\gap{3}\gap{4}\nongap{5}\nongap{6}\nongap{7}\gap{8}\gap{9}\nongap{10} \nongap{11}\nongap{12}\nongap{13}\dots$}  ]
 ]
[.{$\nongap{0}\gap{1}\gap{2}\gap{3}\gap{4}\nongap{5}\nongap{6}\nongap{7}\nongap{8}\gap{9}\nongap{10} \nongap{11}\dots$}  ]
 ]
[.{$\nongap{0}\gap{1}\gap{2}\gap{3}\nongap{4}\gap{5}\generator{6} \generator{7} \nongap{8}\generator{9} \dots$} [.{$\nongap{0}\gap{1}\gap{2}\gap{3}\nongap{4}\gap{5}\gap{6}\generator{7} \nongap{8}\generator{9} \generator{10} \nongap{11}\dots$} [.{$\nongap{0}\gap{1}\gap{2}\gap{3}\nongap{4}\gap{5}\gap{6}\gap{7}\nongap{8} \generator{9} \generator{10} \generator{11} \nongap{12}\nongap{13}\dots$}  ]
[.{$\nongap{0}\gap{1}\gap{2}\gap{3}\nongap{4}\gap{5}\gap{6}\nongap{7}\nongap{8}\gap{9}\generator{10} \nongap{11}\nongap{12}\generator{13} \dots$}  ]
[.{$\nongap{0}\gap{1}\gap{2}\gap{3}\nongap{4}\gap{5}\gap{6}\nongap{7}\nongap{8}\nongap{9}\gap{10}\nongap{11} \nongap{12}\nongap{13}\dots$}  ]
 ]
[.{$\nongap{0}\gap{1}\gap{2}\gap{3}\nongap{4}\gap{5}\nongap{6}\gap{7}\nongap{8} \generator{9} \nongap{10}\generator{11} \dots$} [.{$\nongap{0}\gap{1}\gap{2}\gap{3}\nongap{4}\gap{5}\nongap{6}\gap{7}\nongap{8}\gap{9}\nongap{10} \generator{11} \nongap{12}\generator{13} \dots$}  ]
[.{$\nongap{0}\gap{1}\gap{2}\gap{3}\nongap{4}\gap{5}\nongap{6}\gap{7}\nongap{8}\nongap{9}\nongap{10}\gap{11}\nongap{12} \nongap{13}\dots$}  ]
 ]
[.{$\nongap{0}\gap{1}\gap{2}\gap{3}\nongap{4}\gap{5}\nongap{6}\nongap{7}\nongap{8}\gap{9}\nongap{10} \nongap{11}\dots$}  ]
 ]
[.{$\nongap{0}\gap{1}\gap{2}\gap{3}\nongap{4}\nongap{5}\gap{6}\generator{7} \nongap{8}\nongap{9}\dots$} [.{$\nongap{0}\gap{1}\gap{2}\gap{3}\nongap{4}\nongap{5}\gap{6}\gap{7}\nongap{8} \nongap{9}\nongap{10}\generator{11} \dots$} [.{$\nongap{0}\gap{1}\gap{2}\gap{3}\nongap{4}\nongap{5}\gap{6}\gap{7}\nongap{8}\nongap{9}\nongap{10}\gap{11}\nongap{12} \nongap{13}\dots$}  ]
 ]
 ]
[.{$\nongap{0}\gap{1}\gap{2}\gap{3}\nongap{4}\nongap{5}\nongap{6}\gap{7}\nongap{8} \nongap{9}\dots$}  ]
 ]
[.{$\nongap{0}\gap{1}\gap{2}\nongap{3}\gap{4}\generator{5} \nongap{6}\generator{7} \dots$} [.{$\nongap{0}\gap{1}\gap{2}\nongap{3}\gap{4}\gap{5}\nongap{6} \generator{7} \generator{8} \nongap{9}\dots$} [.{$\nongap{0}\gap{1}\gap{2}\nongap{3}\gap{4}\gap{5}\nongap{6}\gap{7}\generator{8} \nongap{9}\generator{10} \nongap{11}\dots$} [.{$\nongap{0}\gap{1}\gap{2}\nongap{3}\gap{4}\gap{5}\nongap{6}\gap{7}\gap{8}\nongap{9} \generator{10} \generator{11} \nongap{12}\nongap{13}\dots$}  ]
[.{$\nongap{0}\gap{1}\gap{2}\nongap{3}\gap{4}\gap{5}\nongap{6}\gap{7}\nongap{8}\nongap{9}\gap{10}\nongap{11} \nongap{12}\generator{13} \dots$}  ]
 ]
[.{$\nongap{0}\gap{1}\gap{2}\nongap{3}\gap{4}\gap{5}\nongap{6}\nongap{7}\gap{8}\nongap{9} \nongap{10}\generator{11} \dots$} [.{$\nongap{0}\gap{1}\gap{2}\nongap{3}\gap{4}\gap{5}\nongap{6}\nongap{7}\gap{8}\nongap{9}\nongap{10}\gap{11}\nongap{12} \nongap{13}\dots$}  ]
 ]
 ]
[.{$\nongap{0}\gap{1}\gap{2}\nongap{3}\gap{4}\nongap{5}\nongap{6}\gap{7}\nongap{8} \nongap{9}\dots$}  ]
 ]
[.{$\nongap{0}\gap{1}\gap{2}\nongap{3}\nongap{4}\gap{5}\nongap{6} \nongap{7}\dots$}  ]
 ]
[.{$\nongap{0}\gap{1}\nongap{2}\gap{3}\nongap{4} \generator{5} \dots$} [.{$\nongap{0}\gap{1}\nongap{2}\gap{3}\nongap{4}\gap{5}\nongap{6} \generator{7} \dots$} [.{$\nongap{0}\gap{1}\nongap{2}\gap{3}\nongap{4}\gap{5}\nongap{6}\gap{7}\nongap{8} \generator{9} \dots$} [.{$\nongap{0}\gap{1}\nongap{2}\gap{3}\nongap{4}\gap{5}\nongap{6}\gap{7}\nongap{8}\gap{9}\nongap{10} \generator{11} \dots$} [.{$\nongap{0}\gap{1}\nongap{2}\gap{3}\nongap{4}\gap{5}\nongap{6}\gap{7}\nongap{8}\gap{9}\nongap{10}\gap{11}\nongap{12} \generator{13} \dots$}  ]
 ]
 ]
 ]
 ]
 ]
 ] \end{tikzpicture}}

\caption{The tree ${\mathscr T}$ up to depth $7$. White dots refer to the gaps, dark gray dots to the generators and the light gray ones to the elements of the semigroups that are not generators. }\label{tots}
\end{figure}
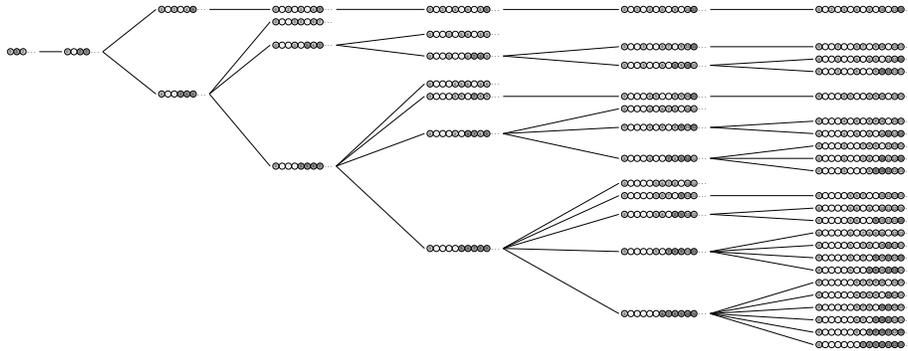

It is well known that all numerical semigroups can be organized in an infinite tree ${\mathscr T}$ whose root is the semigroup ${\mathbb N}_0$ and in which the parent of a numerical semigroup $\Lambda$ is the numerical semigroup $\Lambda'$ obtained by adjoining to $\Lambda$ its Frobenius number. 
For instance,  the parent of the semigroup
$H=\{0,12,19,24,28,31,34,36,38,40,42,43,45,46,47,48,\dots\}$
is the semigroup $H'=\{0,12,19,24,28,31,34,36,38,40,42,43,{\bf 44},45,46,47,48,\dots\}$.
In turn, the children of a numerical semigroup are the semigroups we obtain by taking away one by one the generators that are larger than or equal to the conductor of the semigroup.
The parent of a numerical semigroup of genus $g$ has genus $g-1$ 
and all numerical semigroups are in ${\mathscr T}$, at a depth equal to its genus. In particular, $n_g$ is the number of nodes of ${\mathscr T}$ at depth $g$.
This construction was already considered in 
\cite{RoGaGaJi:oversemigroups}.
Figure~\ref{tots} shows the tree up to depth $7$. 

In \cite{Bras:ordinarization} 
a new tree construction is introduced as follows.
 The {\em ordinarization transform} of a non-ordinary semigroup
 $\Lambda$ with Frobenius number $F$ and multiplicity $m$ is the set
 $\Lambda'=\Lambda\setminus\{m\}\cup\{F\}$.
 For instance,
 the ordinarization transform 
 of the semigroup
$H=\{0,{\bf 12},19,24,28,31,34,36,38,40,42,43,45,46,47,48,\dots\}$
 is the semigroup
 $H'=\{0,19,24,28,31,34,36,38,40,42,43,{\bf 44},45,46,47,48,\dots\}$
The ordinarization transform of an ordinary semigroup is then defined to be itself. 
 Note that the genus of the ordinarization transform of a semigroup is the genus of the semigroup.

The definition of the ordinarization transform of a numerical semigroup allows the construction of a tree 
 ${\mathscr T}_g$ on the set of all numerical semigroups of a given genus rooted at the unique ordinary semigroup of this genus, 
 where the parent of a semigroup is its ordinarization transform and the children of a semigroup are the semigroups obtained by taking away one by one the generators that are larger than the Frobenius number and adding a new non-gap smaller than the multiplicity in a licit place. To illustrate this construction with an example in Figure~\ref{genere} we depicted ${\mathscr T}_7$. 

One significant difference between ${\mathscr T}_g$ and ${\mathscr T}$ is that 
the first one has only a finite number of nodes. In fact, it has $n_g$ nodes, while ${\mathscr T}$ is an infinite tree.
It was conjectured in 
\cite{Bras:ordinarization} 
that the number of numerical semigroups in 
${\mathscr T}_g$ at a given depth 
is at most the number of numerical semigroups in 
${\mathscr T}_{g+1}$ at the same depth.
This was proved in the same reference for the lowest and largest depths.
This conjecture would prove that $n_{g+1}\geq n_{g}$.

\begin{figure}
  \begin{adjustbox}{max totalsize={\textwidth}{\textheight},center}\begin{tikzpicture}[grow=right, every node/.style = {align=left}]\tikzset{level 1+/.style={level distance=12cm}}\Tree[.{$\nongap{0}\gap{1}\gap{2}\gap{3}\gap{4}\gap{5}\gap{6}\gap{7}\generator{8}\generator{9}\generator{10}\generator{11}\generator{12}\generator{13}\generator{14}\dots $} [.{$\nongap{0}\gap{1}\gap{2}\gap{3}\gap{4}\nongap{5}\gap{6}\gap{7}\gap{8}\generator{9}\nongap{10}\generator{11}\generator{12}\generator{13}\nongap{14}\dots $} ]
[.{$\nongap{0}\gap{1}\gap{2}\gap{3}\gap{4}\gap{5}\nongap{6}\gap{7}\gap{8}\generator{9}\generator{10}\generator{11}\nongap{12}\generator{13}\generator{14}\dots $} [.{$\nongap{0}\gap{1}\gap{2}\gap{3}\gap{4}\nongap{5}\nongap{6}\gap{7}\gap{8}\gap{9}\nongap{10}\nongap{11}\nongap{12}\generator{13}\generator{14}\dots $} ]
[.{$\nongap{0}\gap{1}\gap{2}\nongap{3}\gap{4}\gap{5}\nongap{6}\gap{7}\gap{8}\nongap{9}\gap{10}\generator{11}\nongap{12}\generator{13}\nongap{14}\dots $} ]
[.{$\nongap{0}\gap{1}\gap{2}\nongap{3}\gap{4}\gap{5}\nongap{6}\gap{7}\gap{8}\nongap{9}\nongap{10}\gap{11}\nongap{12}\nongap{13}\generator{14}\dots $} ]
[.{$\nongap{0}\gap{1}\gap{2}\gap{3}\gap{4}\nongap{5}\nongap{6}\gap{7}\gap{8}\nongap{9}\nongap{10}\nongap{11}\nongap{12}\gap{13}\nongap{14}\dots $} ]
]
[.{$\nongap{0}\gap{1}\gap{2}\gap{3}\gap{4}\gap{5}\gap{6}\nongap{7}\gap{8}\generator{9}\generator{10}\generator{11}\generator{12}\generator{13}\nongap{14}\dots $} [.{$\nongap{0}\gap{1}\gap{2}\gap{3}\gap{4}\nongap{5}\gap{6}\nongap{7}\gap{8}\gap{9}\nongap{10}\generator{11}\nongap{12}\generator{13}\nongap{14}\dots $} ]
[.{$\nongap{0}\gap{1}\gap{2}\gap{3}\gap{4}\gap{5}\nongap{6}\nongap{7}\gap{8}\gap{9}\generator{10}\generator{11}\nongap{12}\nongap{13}\nongap{14}\dots $} ]
[.{$\nongap{0}\gap{1}\gap{2}\gap{3}\gap{4}\gap{5}\nongap{6}\nongap{7}\gap{8}\nongap{9}\gap{10}\generator{11}\nongap{12}\nongap{13}\nongap{14}\dots $} ]
[.{$\nongap{0}\gap{1}\gap{2}\gap{3}\gap{4}\nongap{5}\gap{6}\nongap{7}\gap{8}\nongap{9}\nongap{10}\gap{11}\nongap{12}\generator{13}\nongap{14}\dots $} ]
[.{$\nongap{0}\gap{1}\gap{2}\gap{3}\gap{4}\gap{5}\nongap{6}\nongap{7}\gap{8}\nongap{9}\nongap{10}\gap{11}\nongap{12}\nongap{13}\nongap{14}\dots $} ]
[.{$\nongap{0}\gap{1}\gap{2}\gap{3}\gap{4}\nongap{5}\gap{6}\nongap{7}\gap{8}\nongap{9}\nongap{10}\nongap{11}\nongap{12}\gap{13}\nongap{14}\dots $} ]
]
[.{$\nongap{0}\gap{1}\gap{2}\gap{3}\nongap{4}\gap{5}\gap{6}\gap{7}\nongap{8}\gap{9}\generator{10}\generator{11}\nongap{12}\generator{13}\nongap{14}\dots $} ]
[.{$\nongap{0}\gap{1}\gap{2}\gap{3}\gap{4}\nongap{5}\gap{6}\gap{7}\nongap{8}\gap{9}\nongap{10}\generator{11}\generator{12}\nongap{13}\generator{14}\dots $} ]
[.{$\nongap{0}\gap{1}\gap{2}\gap{3}\gap{4}\gap{5}\nongap{6}\gap{7}\nongap{8}\gap{9}\generator{10}\generator{11}\nongap{12}\generator{13}\nongap{14}\dots $} [.{$\nongap{0}\gap{1}\gap{2}\gap{3}\nongap{4}\gap{5}\nongap{6}\gap{7}\nongap{8}\gap{9}\nongap{10}\gap{11}\nongap{12}\generator{13}\nongap{14}\dots $} [.{$\nongap{0}\gap{1}\nongap{2}\gap{3}\nongap{4}\gap{5}\nongap{6}\gap{7}\nongap{8}\gap{9}\nongap{10}\gap{11}\nongap{12}\gap{13}\nongap{14}\dots $} ]
]
[.{$\nongap{0}\gap{1}\gap{2}\gap{3}\nongap{4}\gap{5}\nongap{6}\gap{7}\nongap{8}\gap{9}\nongap{10}\nongap{11}\nongap{12}\gap{13}\nongap{14}\dots $} ]
]
[.{$\nongap{0}\gap{1}\gap{2}\gap{3}\gap{4}\gap{5}\gap{6}\nongap{7}\nongap{8}\gap{9}\generator{10}\generator{11}\generator{12}\generator{13}\nongap{14}\dots $} [.{$\nongap{0}\gap{1}\gap{2}\gap{3}\nongap{4}\gap{5}\gap{6}\nongap{7}\nongap{8}\gap{9}\gap{10}\nongap{11}\nongap{12}\generator{13}\nongap{14}\dots $} ]
[.{$\nongap{0}\gap{1}\gap{2}\gap{3}\gap{4}\gap{5}\nongap{6}\nongap{7}\nongap{8}\gap{9}\gap{10}\generator{11}\nongap{12}\nongap{13}\nongap{14}\dots $} ]
[.{$\nongap{0}\gap{1}\gap{2}\gap{3}\gap{4}\nongap{5}\gap{6}\nongap{7}\nongap{8}\gap{9}\nongap{10}\gap{11}\nongap{12}\nongap{13}\nongap{14}\dots $} ]
[.{$\nongap{0}\gap{1}\gap{2}\gap{3}\gap{4}\gap{5}\nongap{6}\nongap{7}\nongap{8}\gap{9}\nongap{10}\gap{11}\nongap{12}\nongap{13}\nongap{14}\dots $} ]
[.{$\nongap{0}\gap{1}\gap{2}\gap{3}\nongap{4}\gap{5}\gap{6}\nongap{7}\nongap{8}\gap{9}\nongap{10}\nongap{11}\nongap{12}\gap{13}\nongap{14}\dots $} ]
]
[.{$\nongap{0}\gap{1}\gap{2}\gap{3}\nongap{4}\gap{5}\gap{6}\gap{7}\nongap{8}\nongap{9}\gap{10}\generator{11}\nongap{12}\nongap{13}\generator{14}\dots $} ]
[.{$\nongap{0}\gap{1}\gap{2}\gap{3}\gap{4}\gap{5}\nongap{6}\gap{7}\nongap{8}\nongap{9}\gap{10}\generator{11}\nongap{12}\generator{13}\nongap{14}\dots $} [.{$\nongap{0}\gap{1}\gap{2}\nongap{3}\gap{4}\gap{5}\nongap{6}\gap{7}\nongap{8}\nongap{9}\gap{10}\nongap{11}\nongap{12}\gap{13}\nongap{14}\dots $} ]
]
[.{$\nongap{0}\gap{1}\gap{2}\gap{3}\gap{4}\gap{5}\gap{6}\nongap{7}\nongap{8}\nongap{9}\gap{10}\generator{11}\generator{12}\generator{13}\nongap{14}\dots $} [.{$\nongap{0}\gap{1}\gap{2}\gap{3}\gap{4}\gap{5}\nongap{6}\nongap{7}\nongap{8}\nongap{9}\gap{10}\gap{11}\nongap{12}\nongap{13}\nongap{14}\dots $} ]
]
[.{$\nongap{0}\gap{1}\gap{2}\gap{3}\nongap{4}\gap{5}\gap{6}\gap{7}\nongap{8}\nongap{9}\nongap{10}\gap{11}\nongap{12}\nongap{13}\nongap{14}\dots $} ]
[.{$\nongap{0}\gap{1}\gap{2}\gap{3}\gap{4}\nongap{5}\gap{6}\gap{7}\nongap{8}\nongap{9}\nongap{10}\gap{11}\generator{12}\nongap{13}\nongap{14}\dots $} ]
[.{$\nongap{0}\gap{1}\gap{2}\gap{3}\gap{4}\gap{5}\nongap{6}\gap{7}\nongap{8}\nongap{9}\nongap{10}\gap{11}\nongap{12}\generator{13}\nongap{14}\dots $} ]
[.{$\nongap{0}\gap{1}\gap{2}\gap{3}\gap{4}\gap{5}\gap{6}\nongap{7}\nongap{8}\nongap{9}\nongap{10}\gap{11}\generator{12}\generator{13}\nongap{14}\dots $} ]
[.{$\nongap{0}\gap{1}\gap{2}\gap{3}\gap{4}\nongap{5}\gap{6}\gap{7}\nongap{8}\nongap{9}\nongap{10}\nongap{11}\gap{12}\nongap{13}\nongap{14}\dots $} ]
[.{$\nongap{0}\gap{1}\gap{2}\gap{3}\gap{4}\gap{5}\gap{6}\nongap{7}\nongap{8}\nongap{9}\nongap{10}\nongap{11}\gap{12}\generator{13}\nongap{14}\dots $} ]
[.{$\nongap{0}\gap{1}\gap{2}\gap{3}\gap{4}\gap{5}\nongap{6}\gap{7}\nongap{8}\nongap{9}\nongap{10}\nongap{11}\nongap{12}\gap{13}\nongap{14}\dots $} ]
[.{$\nongap{0}\gap{1}\gap{2}\gap{3}\gap{4}\gap{5}\gap{6}\nongap{7}\nongap{8}\nongap{9}\nongap{10}\nongap{11}\nongap{12}\gap{13}\nongap{14}\dots $} ]
]
 \end{tikzpicture}\end{adjustbox}

\caption{The whole tree ${\mathscr T}_7$}\label{genere}
\end{figure}

In Section~\ref{s:q} we will construct the quasi-ordinarization transform of a general semigroup, paralleling the ordinarization transform. If the quasi-ordinarization transform is applied repeatedly to a numerical semigroup, it ends up in a quasi-ordinary semigroup.
In Section~\ref{s:qn} we define the quasi-ordinarization number of a semigroup as the number of successive quasi-ordinarization transforms of the semigroup that give a quasi-ordinary semigroup.
Section~\ref{s:varrho} analyzes the number of numerical semigroups of a given genus and a given quasi-ordinarization number in terms of the given parameters.
We present the conjecture that the number of numerical semigroups of a given genus and a fixed quasi-ordinarization number increases with the genus and we prove it for the largest quasi-ordinarization numbers.
In Section~\ref{s:forest} we present the forest of semigroups of a given genus that is obtained when connecting each semigroup to its quasi-ordinarization transform. The forest corresponding to genus $g$ is denoted ${\mathscr F}_g$.
Section~\ref{s:Ts} analyzes the relationships between ${\mathscr T}$, ${\mathscr T}_g$, and ${\mathscr F}_g$. 

From the perspective of the forests of numerical semigroups here presented, the conjecture in Section~\ref{s:varrho} translates to the conjecture that the number of numerical semigroups in 
${\mathscr F}_g$ at a given depth 
is at most the number of numerical semigroups in 
${\mathscr F}_{g+1}$ at the same depth.
The results in Section~\ref{s:varrho} provide a proof of the conjecture for the largest depths. Proving this conjecture for all depths, would prove that $n_{g+1}\geq n_{g}$.
Hence, we expect our work to contribute to the proof of the conjectured increasingness of the sequence $n_g$ (\href{http://oeis.org/A007323}{A007323}).

\section{Quasi-ordinary semigroups and quasi-ordinarization transform}
\label{s:q}

  {\em Quasi-ordinary} semigroups are those semigroups for which $m=g$ and so, there is a unique gap larger than $m$.
  The {\em sub-Frobenius number} of a non-ordinary semigroup $\Lambda$ with Frobenius number $F$ is the Frobenius number of $\Lambda\cup\{F\}$.

  The {\em subconductor} of a semigroup with Frobenius number $F$ is the smallest nongap in the interval of nongaps immediatelly previous to $F$. For instance, the subconductor of the above example,
$H=\{0,12,19,24,28,31,34,36,38,40,42,43,45,46,47,48,\dots\}$, is $42$.

\begin{Lemma}
  Let $\Lambda$ be a non-ordinary and non quasi-ordinary semigroup, with multiplicity $m$, genus $g$, and sub-Frobenius number $f$. Then
  $\Lambda\cup\{f\}\setminus \{m\}$ is another numerical semigroup of the same genus $g$.
\end{Lemma}

\begin{proof}
Since $\Lambda$ is already a numerical semigroup, it is enough to see that
$F-f$ is not in $\Lambda\cup\{f\}\setminus \{m\}$, where $F$ is
  the Frobenius number of $\Lambda$.
  Notice that for a non-ordinary numerical semigroup, the difference between its Frobenius number and its sub-Frobenius number needs to be less than the multiplicity of the semigroup; hence, $F-f\not\in \Lambda$. So the only option for $F-f$ to be in $\Lambda\cup\{f\}\setminus \{m\}$ is that $F-f=f$. In this case, any integer between $1$ and $f-1$ must be a gap, since the integers between $F-1$ and $F-f+1$ are nongaps. In this case, $\Lambda$ would be quasi-ordinary, contradicting the hypotheses.
  \end{proof}

\begin{Definition}
The {\em quasi-ordinarization transform} of a non-ordinary and non quasi-ordinary numerical semigroup $\Lambda$, with
multiplicity $m$, genus $g$ and sub-Frobenius number $f$,
is the numerical semigroup $\Lambda\cup\{f\}\setminus \{m\}$.

The {\em quasi-ordinarization} of either an ordinary or quasi-ordinary
semigroup is defined to be itself.
\end{Definition}

As an example, the quasi-ordinarization of the well-tempered harmonic semigroup 
$H=\{0,{\bf 12},19,24,28,31,34,36,38,40,42,43,45,46,47,48,\dots\}$ used in the previous examples 
is $H'=\{0,19,24,28,31,34,36,38,40,{\bf 41},42,43,45,46,47,48,\dots\}.$

\begin{Remark} In the ordinarization and quasi-ordinarization transform process, we replace the multiplicity by the largest and second largest gap respectively, and we obtain numerical semigroups.
In general, if we replace the multiplicity by the third largest gap,
we do not obtain a numerical semigroup.

See for instance $\{0,2,4,6,8,10,11,\dots\}$.
Replacing $2$ by $5$, we obtain
$\{0,4,5,6,8,10,11,\dots\}$ which is not a numerical semigroup since 
$9=4+5$ is not in the set.
\end{Remark}

\section{Quasi-ordinarization number}
\label{s:qn}

Next lemma explicits that there is only one quasi-ordinary semigroup with genus 
$g$ and conductor $c$ where $c\leq 2g$.

\begin{Lemma}\label{l:roots}
For each positive integers $g$ and $c$ with $c\leq 2g$ 
the semigroup $\{0,g,g+1,\dots, c-2,c,c+1\dots \}$
is the unique quasi-ordinary semigroup of genus $g$ and conductor $c$.
\end{Lemma}

  The quasi-ordinarization transform of a non-ordinary semigroup
  of genus $g$ and conductor $c$ can be applied subsequently and
  at some step we will attain the quasi-ordinary semigroup of that genus and conductor, that is, the numerical semigroup
  $\{0,g,g+1,\dots,c-2,c,c+1,\dots\}$.
  The number of such steps is defined to be the {\em quasi-ordinarization number} of $\Lambda$. 

  We denote by $\varrho_{g,q}$ the number of numerical semigroups of genus $g$ and quasi-ordinarization number $q$.
  In Table~\ref{t:ngq} one can see the values of $\varrho_{g,q}$ for genus up to 45. It has been computed by an exhaustive exploration of the semigroup tree using the RGD algorithm \cite{rgd}.
  
\begin{Lemma}
\label{lemma:equivdef}
The quasi-ordinarization number of a non-ordinary numerical semigroup of genus $g$ coincides with the number of non-zero non-gaps of the semigroup that are smaller than or equal to $g-1$. 
\end{Lemma}

\begin{proof}
A non-ordinary numerical semigroup of genus $g$ 
is non-quasi-ordinary if and only if its multiplicity 
is at most $g-1$. 
Consequently, we can repeatedly apply the quasi-ordinarization transform to a numerical semigroup while its multiplicity is at most $g-1$. Furthermore, the number of consecutive transforms that we can apply before obtaining the quasi-ordinary semigroup is hence the number of its non-zero non-gaps that are at most the genus minus one. 
\end{proof}

For a numerical semigroup $\Lambda$ we will consider its enumeration $\lambda$, that is, the unique increasing bijective map between ${\mathbb N}_0$ and $\Lambda$. The element $\lambda(i)$ is then denoted $\lambda_i$.
As a consequence of the previous lemma, for a numerical semigroup $\Lambda$ with quasi-ordinarization number equal to $q$, the non-gaps that are at most $g-1$ are exactly $\lambda_0=0, \lambda_1, \dots, \lambda_q$.

\begin{Lemma}
The maximum quasi-ordinarization number of a non-ordinary semigroup of genus $g$ is~$\lfloor\frac{g-1}{2}\rfloor$. 
\end{Lemma}

\begin{proof}
Let $\Lambda$ be a numerical semigroup with quasi-ordinarization number equal to $q$. 
Since the Frobenius number $F$ is at most $2g-1$,
the total number of gaps 
from $1$ to $2g-1$ is $g$
and so the number of non-gaps  
from $1$ to $2g-1$ is $g-1$. 
The number of those non-gaps 
that are larger than $g-1$ is $g-1-q$. 
On the other hand 
$\lambda_q+\lambda_1, \lambda_q+\lambda_2, \dots, 2\lambda_q$
are different non-gaps between $g$ and $2g-1$. 
So the number of non-gaps between $g$ and $2g-1$ is at least $q$. 
All these results imply that $g-1-q\geq q$ and so, $q\leq\frac{g-1}{2}$. 

On the other hand, the bound stated in the lemma is attained by the hyperelliptic numerical semigroup 
\begin{equation}
\{0,2,4,\dots,2\left\lfloor\frac{g-1}{2}\right\rfloor,2\left(\left\lfloor\frac{g-1}{2}\right\rfloor +1\right),\dots,2g,2g+1,2g+2,\dots\}. 
\label{hyperelliptic}
\end{equation}
\end{proof}

We will next see that the maximum ordinarization number stated in the previous lemma is attained uniquely by 
the numerical semigroup in (\ref{hyperelliptic}). 
To prove this result we will need the next lemma.
Let us recall that $A+B=\{a+b : a\in A, b\in B\}$ and that $\#A$  denotes the cardinality of $A$. 

 \begin{Lemma}
 \label{lemma:arithseqsums}
 Consider a finite subset $A=\{a_1<\dots<a_n\}\subseteq {\mathbb N}_0$. 
 \begin{enumerate}
 \item 
 The set $A+A$ contains at least $2n-1$ elements 
 \item If $n\geq 1$, 
   the set $A+A$ contains exactly $2n-1$ elements if and only if
   there exists a positive integer $\alpha$  such that
 $a_{i}=a_1+(i-1)\alpha$ for all $i\leq n$.
 \item If $n\geq 4$, the set $A+A$ contains exactly $2n$ elements if and only if either
   \begin{itemize}
     \item there exists a positive integer $\alpha$ such that
       $a_{i}=a_1+\alpha (i-1)$ for all $i$ with $1\leq i<n$ and $a_n=a_1+n\alpha$,
     \item there exists a positive integer $\alpha$ such that
       $a_{i}=a_1+i\alpha$ for all $i$ with $2\leq i\leq n$.
       \end{itemize}
 \end{enumerate}
 \end{Lemma}

 \begin{proof}

   The first item stems from the fact that if
   $A=\{a_1,\dots,a_n\}$, then $A+A$ must contain at least
   $2a_1,a_1+a_2,a_1+a_3,\dots,a_1+a_n,a_2+a_n,a_3+a_n,\dots,a_{n-1}+a_n,2a_n$, which are all different.

   The second item easily follows from the fact that if $A+A$ has $2n-1$ elements, then $A+A$ must be exactly the set
   $2a_1,a_1+a_2,a_1+a_3,\dots,a_1+a_n,a_2+a_n,a_3+a_n,\dots,a_{n-1}+a_n,2a_n$.
   Indeed, in this case the increasing set $\{a_1+a_3,\dots,a_1+a_n,a_2+a_n,a_3+a_n,\dots,a_{n-1}+a_n,2a_n\}$ must coincide with the increasing set
   $\{2a_2,a_2+a_3,a_2+a_4,\dots,a_2+a_n,a_3+a_n,\dots,a_{n-1}+a_n,2a_n\}$, having as a consequence that $2a_2=a_1+a_3$ and so, $a_2=\frac{a_1+a_3}{2}=a_1+\frac{a_3-a_1}{2}$, and $a_3=2a_2-a_1=a_1+2\frac{a_3-a_1}{2}$.
   Hence,
   \begin{eqnarray*}
     a_2&=&a_1+\frac{a_3-a_1}{2}\\
     a_3&=&a_1+2\frac{a_3-a_1}{2}\\
   \end{eqnarray*}
   Similarly one can show that $2a_3=a_2+a_4$ and, so, $a_4=2a_3-a_2=2a_1+4\frac{a_3-a_1}{2}-a_1-\frac{a_3-a_1}{2}=a_1+3\frac{a_3-a_1}{2}$.
   And it equally follows that
   \begin{eqnarray*}
     a_4&=&a_1+3\frac{a_3-a_1}{2}\\
     a_5&=&a_1+4\frac{a_3-a_1}{2}\\
     &\vdots&\\
   \end{eqnarray*}

For the third item, one direction of the proof is obvious, so  we just need to prove the other one, that is, if the sum contains $2n$ elements, then $a_1,\dots,a_n$ must be as stated.
   
   We will proceed by induction. Suppose that $n=4$ and that the set $A+A$ contains exactly $8$ elements. Since the ordered sequence
   \begin{equation}\label{e:szero}2a_1<a_1+a_2<2a_2<a_2+a_3<2a_3<a_3+a_4<2a_4\end{equation}
   contains already $7$ elements, then necessarily two of the elements 
$a_1+a_3, a_1+a_4, a_2+a_4$ coincide with one element in \eqref{e:szero} and the third one is not in \eqref{e:szero}.
   So, at least one of $a_1+a_3$ and $a_2+a_4$ must be in \eqref{e:szero}.

   Suppose first that $a_1+a_3$ is in \eqref{e:szero}.
   Then necessarily $a_1+a_3=2a_2$ which means that $a_2-a_1=a_3-a_2$. Hence,
   there exists $\alpha$ (in fact, $\alpha=a_2-a_1$)
   such that $a_2=a_1+\alpha$ and $a_3=a_1+2\alpha$.
   Now, the elements 
   \begin{equation}\label{e:sone}2a_1<a_1+a_2<2a_2<a_2+a_3<2a_3\end{equation}
   are equally separated by the same separation $\alpha$.
   That is,
   \begin{eqnarray*}
     (a_1+a_2)-(2a_1)&=&\alpha\\
     (2a_2)-(a_1+a_2)&=&\alpha\\
     (a_2+a_3)-(2a_2)&=&\alpha\\
     (2a_3)-(a_2+a_3)&=&\alpha.\\
   \end{eqnarray*}
   And also the elements
   \begin{equation}\label{e:stwo}a_4+a_1<a_4+a_2<a_4+a_3\end{equation}
   are equally separated by the same separation $\alpha$.
   That is,
   \begin{eqnarray*}
     (a_4+a_3)-(a_4+a_2)&=&\alpha\\
     (a_4+a_2)-(a_4+a_1)&=&\alpha.\\
   \end{eqnarray*}

   Furthermore, $A+A$ must contain all the elements in \eqref{e:sone} and \eqref{e:stwo} as well as the element $2a_4$, which is not in \eqref{e:sone}, nor in \eqref{e:stwo}. Since $\#(A+A)=8$, this means that there must be exatly one element that is both in \eqref{e:sone} and in \eqref{e:stwo}. The only way for this to happen is that $2a_3=a_4+a_1$. Consequently, $a_4+a_1=2a_1+4\alpha$, and so, $a_4=a_1+4\alpha$. This proves the result in the first case.

   For the case in which $a_2+a_4$ is in \eqref{e:szero}, necessarily $a_2+a_4=2a_3$, which means that $a_3-a_2=a_4-a_3$.
   Hence,
   there exists $\beta$ (in fact, $\beta=a_3-a_2$)
   such that $a_3=a_2+\beta$ and $a_4=a_2+2\beta$.
      Now, the elements 
   \begin{equation}\label{e:sonebeta}2a_2<a_2+a_3<2a_3<a_3+a_4<2a_4\end{equation}
   are equally separated by the same separation $\beta$.
   That is,
   \begin{eqnarray*}
     (a_2+a_3)-(2a_2)&=&\beta\\
     (2a_3)-(a_2+a_3)&=&\beta\\
     (a_3+a_4)-(2a_3)&=&\beta\\
     (2a_4)-(a_3+a_4)&=&\beta.\\
   \end{eqnarray*}
   And also the elements
   \begin{equation}\label{e:stwobeta}   a_1+a_2<a_1+a_3<a_1+a_4\end{equation}
   are equally separated by the same separation $\beta$.
   That is,
   \begin{eqnarray*}
     (a_1+a_3)-(a_1+a_2)&=&\beta\\
     (a_1+a_4)-(a_1+a_3)&=&\beta.\\
   \end{eqnarray*}
   Now, $A+A$ must contain all the elements in \eqref{e:sonebeta} and \eqref{e:stwobeta} as well as the element $2a_1$, which is not in \eqref{e:sonebeta}, nor in \eqref{e:stwobeta}. Since $\#(A+A)=8$, this means that there must be exactly one element that is both in \eqref{e:sonebeta} and in \eqref{e:stwobeta}. The only way for this to happen is that $a_1+a_4=2a_2$.
   Consequently, $a_1+a_4=2a_1+4\alpha$, and so, $a_4=a_1+4\alpha$.
Hence, $a_2=a_1+2\beta$, $a_3=a_1+3\beta$, $a_4=a_1+4\beta$.
   This proves the result in the second case and concludes the proof for $n=4$.

   Now, let us prove the result for a general $n>4$.
   We will denote $A_n$ the set $\{a_1,\dots,a_n\}$.

   Notice that $A_1+A_1=\{2a_1\}$ while, if $i>1$, then $\{a_{i-1}+a_i,2a_i\}\subseteq (A_i+A_i)\setminus(A_{i-1}+A_{i-1})$, hence, $\#((A_i+A_i)\setminus(A_{i-1}+A_{i-1}))\geq 2$.
   Consequently, if $\#(A_n+A_n)=2n,$ we can afirm that there exists exactly one integer $s$ such that
$\#(A_{r}+A_{r})=2r-1,$ for all $r<s$ and  
$\#(A_{r}+A_{r})=2r$ for all $r\geq s$.

If $s=n$, then we already have, by the second item of the lemma, that  $a_{i}=a_1+(i-1)\gamma$ for a given positive integer $\gamma$ for all $i<n$.
   
   On one hand, 
   \begin{equation}\label{e:setone}A_{n-1}+A_{n-1}=\{2a_1,2a_1+\gamma,2a_1+2\gamma,2a_1+3\gamma,\dots,2a_1+(2n-4)\gamma\},\end{equation}
   which has $2n-3$ elements.
   On the other hand,
   \begin{equation}\label{e:settwo}A_{n-1}+a_{n}=\{(a_1+a_n),(a_1+a_n)+\gamma,(a_1+a_n)+2\gamma,(a_1+a_n)+3\gamma,\dots,(a_1+a_n)+(n-2)\gamma\},\end{equation}
has $n-1$ elements.
   
Now, $A+A=(A_{n-1}+A_{n-1})\cup (A_{n-1}+a_n)\cup (2a_n)$.
By the inclusion-exclusion principle, and since $2a_n$ is not in $(A_{n-1}+A_{n-1})\cup (A_{n-1}+a_n)$,

\begin{eqnarray*}\#\left((A_{n-1}+A_{n-1})\cap (A_{n-1}+a_n)\right)&=&\#(A_{n-1}+A_{n-1})+\#(A_{n-1}+a_n)+1-\#(A+A)\\
  &=& (2n-3)+(n-1)+1-2n\\
  &=&n-3
  \end{eqnarray*}

By \eqref{e:setone} and \eqref{e:settwo},
we conclude that $(a_1+a_n)+(n-4)\gamma=2a_1+(2n-4)\gamma$, that is,
$a_n=a_1+n\gamma$.
Hence the result follows with $\alpha=\gamma$.

On the contrary, if $s<n$, then, since $\#\left(A_{n-1}+A_{n-1}\right)=2(n-1)$, we can apply the induction hypothesis and affirm that
either one of  the following cases, (a) or (b), holds.
   \begin{enumerate}
     \item[(a)] There exists a positive integer $\alpha_{n-1}$ such that
       $a_{i}=a_1+\alpha_{n-1} (i-1)$ for all $i$ with $1\leq i<n-1$ and $a_{n-1}=a_1+(n-1)\alpha$;
     \item[(b)] There exists a positive integer $\alpha$ such that
       $a_{i}=a_1+i\alpha_{n-1}$ for all $i$ with $2\leq i\leq n-1$.
   \end{enumerate}

In case (a) we will have
$$
    \begin{array}{r}
          A_{n-1}+A_{n-1}=
          \{2a_1,2a_1+\alpha_{n-1},2a_1+2\alpha_{n-1},\dots\\\dots,2a_1+(2n-4)\alpha_{n-1}, 2a_1+(2n-2)\alpha_{n-1}\},
          \end{array}
$$
    and
    $$
  \begin{array}{r}
    A_{n-1}+a_{n}=
        \{(a_1+a_n),(a_1+a_n)+\alpha_{n-1},(a_1+a_n)+2\alpha_{n-1},\dots\\\dots,(a_1+a_n)+(n-3)\alpha_{n-1},(a_1+a_n)+(n-1)\alpha_{n-1}\},
\end{array}
$$

In case (b) we will have
$$
  \begin{array}{r}
          A_{n-1}+A_{n-1}=
          \{2a_1,2a_1+2\alpha_{n-1},2a_1+3\alpha_{n-1},\dots\\\dots,2a_1+(2n-3)\alpha_{n-1},2a_1+(2n-2)\alpha_{n-1}\},
          \end{array}
$$
and
$$
    \begin{array}{r}
  A_{n-1}+a_{n}=
  \{(a_1+a_n),(a_1+a_n)+2\alpha_{n-1},(a_1+a_n)+3\alpha_{n-1},\dots\\\dots,(a_1+a_n)+(n-1)\alpha_{n-1}\},
    \end{array}
    $$

Now,  
\begin{eqnarray*}\#\left((A_{n-1}+A_{n-1})\cap (A_{n-1}+a_n)\right)&=&\#(A_{n-1}+A_{n-1})+\#(A_{n-1}+a_n)+1-\#(A+A)\\
  &=& \#(A_{n-1}+A_{n-1})-n\\
  &=& n-2.
\end{eqnarray*}

This is only possible in case (b) with

$$
    \begin{array}{r}
(A_{n-1}+A_{n-1})\cap (A_{n-1}+a_n)=
      \{(a_1+a_2),(a_1+a_n)+2\alpha_{n-1},(a_1+a_n)+3\alpha_{n-1},\dots\\\dots,(a_1+a_n)+(n-2)\alpha_{n-1}\},
    \end{array}
$$

and, hence, with
$(a_1+a_n)+(n-2)\alpha_{n-1}=2a_1+(2n-2)\alpha_{n-1}$, that is,
$a_n=a_1+n\alpha_{n-1}$, hence yielding the result with $\alpha=\alpha_{n-1}$.

 \end{proof}

 \begin{Lemma}
 Let $g>2$ and $g\neq 4, g\neq 6$. 
 The unique non-quasi-ordinary numerical semigroup of genus $g$ and 
 quasi-ordinarization number $\lfloor\frac{g-1}{2}\rfloor$ is $\{0,2,4,\dots,2g,2g+2,2g+3\dots\}$. 
 \end{Lemma}

\begin{proof}
  If $g=3$, there is only one numerical semigroup non-ordinary and non-quasi-ordinary as we can observe in Figure~\ref{tots}, and it is exactly $\{0,2,4,6,\dots\}$, which indeed, has quasi-ordinarization number $\lfloor\frac{g-1}{2}\rfloor$ and it is of the form $\{0,2,4,\dots,2g,2g+1,2g+2,\dots\}$.
The case $g=4$ and $g=6$ are excluded from the statement (and analyzed in Remark~\ref{r:casquatre}).
So, we can assume that either $g=5$ or $g>6$.

Suppose that the quasi-ordinarization number of $\Lambda$ is $\lfloor\frac{g-1}{2}\rfloor$. 
Since $\lambda_{\lfloor\frac{g-1}{2}\rfloor}\leq g-1$,
we know that the set of all non-gaps between $0$ and $2g-2$ 
must contain all the sums 
$$\Sigma=\{\lambda_i+\lambda_j:0\leq i,j\leq \lfloor\frac{g-1}{2}\rfloor\}.$$
But the number of non-gaps between $0$ and $2g-2$ is either $g-1$ or
$g$ depending on whether $2g-1$ is a gap or not.
So, $\#\Sigma\leq g$.
On the other hand, by
Lemma~\ref{lemma:arithseqsums},
$\#\Sigma\geq 2\lfloor\frac{g-1}{2}\rfloor+1$. 

If $g$ is odd, we get that 
$2\lfloor\frac{g-1}{2}\rfloor+1=g$ and so, $\#\Sigma=g$. Then, by the second item in Lemma~\ref{lemma:arithseqsums},
we get that 
$\lambda_i=i\lambda_1$ 
for $i\leq\frac{g-1}{2}$. 
Now $\lambda_{\frac{g-1}{2}}=\frac{g-1}{2}\lambda_1$ and, since $\lambda_{\frac{g-1}{2}}\leq g-1$, one can deduce that
$\lambda_1\leq 2$. If $\lambda_1=1$ this contradicts 
$g>1$. 
So,
$\lambda_i=2i$ 
for $0\leq i\leq\frac{g-1}{2}$ and the 
remaining non-gaps 
between $g$ and $2g$ are necessarily 
$\lambda_i=2i$ 
for $i=\frac{g-1}{2}+1$ to $i=g$. 

If $g$ is even then $g-1\leq\#\Sigma\leq g$.
If $\#\Sigma=g$, then, since the number of summands in the sum $\Sigma$ is at least $4$ (because we excluded the even genera $4$ and $6$), we can apply
the third item in
Lemma~\ref{lemma:arithseqsums}. Then we get that
$\lambda_{\frac{g}{2}-1}= \frac{g}{2}\lambda_1$. This, together with
$\lambda_{\frac{g}{2}-1}\leq g-1$ implies that
$\lambda_1\leq 2\frac{g-1}{g}<2$. So, $\lambda_1=1$, contradicting
$g>1$.
Hence, it must be $\Sigma=g-1$.
If $\Sigma=g-1$, then,
by the second item in
Lemma~\ref{lemma:arithseqsums},
we get that 
$\lambda_i=i\lambda_1$ 
for all $i\leq\frac{g}{2}-1$.
Now $\lambda_{\frac{g}{2}-1}= 
(\frac{g}{2}-1)\lambda_1$
and, since $\lambda_{\frac{g}{2}-1}\leq g-1$, one can deduce that
$\lambda_1\leq 2\frac{g-1}{g-2}$.
But $2\frac{g-1}{g-2}<3$ if $g\geq 5$.
So, $\lambda_1$ con only be $1$ or $2$.
If $\lambda_1=1$ this contradicts 
$g>1$. 
So,
$\lambda_i=2i$ 
for $0\leq i\leq\frac{g}{2}-1$ and the 
remaining non-gaps 
between $g$ and $2g$ are necessarily 
$\lambda_i=2i$ 
for $i=\frac{g}{2}$ to $i=g$.

\end{proof}

\begin{Remark}
  \label{r:casquatre}
  For $g=4$, the maximum quasi-ordinarization number $\lfloor\frac{g-1}{2}\rfloor=1$ is, in fact, attained by three of the 7 semigroups of genus 4.
The semigroups whose quasi-ordinarization number is maximum are 
$\{0,3,6,\dots\}$,
$\{0,2,4,6,8,\dots\}$,
$\{0,3,5,6,8,\dots\}$.

  For $g=6$, the maximum quasi-ordinarization number $\lfloor\frac{g-1}{2}\rfloor=2$ is, in fact, attained by two of the 23 semigroups of genus 6.
The semigroups whose quasi-ordinarization number is maximum are 
$\{0,2,4,6,8,10,12,\dots\}$,
$\{0,4,5,8,9,10,12,\dots\}$.

Hence, $g=4$ and $g=6$ are exceptional cases.
\end{Remark}

{
  \begin{table}
    {  \tiny   \begin{center}
\resizebox{\textwidth}{!}{\begin{tabular}{c}
\begin{tabular}{|c|c|c|c|c|c|c|c|c|c|c|}
\hline $\varrho_{g,q}$&$g=1$ &$g=2$ &$g=3$ &$g=4$ &$g=5$ &$g=6$ &$g=7$ &$g=8$ &$g=9$ &$g=10$ \\\hline 
$q=0$ &1 &2 &3 &4 &5 &6 &7 &8 &9 &10 \\ 
$q=1$ &&&1 &3 &6 &15 &24 &42 &61 &93 \\ 
$q=2$ &&&&&1 &2 &7 &16 &43 &89 \\ 
$q=3$ &&&&&&&1 &1 &4 &11 \\ 
$q=4$ &&&&&&&&&1 &1 \\ 
\hline
sum= &1 &2 &4 &7 &12 &23 &39 &67 &118 &204 \\\hline
\end{tabular}\\

\begin{tabular}{|c|c|c|c|c|c|c|c|c|c|c|}
\hline $\varrho_{g,q}$&$g=11$ &$g=12$ &$g=13$ &$g=14$ &$g=15$ &$g=16$ &$g=17$ &$g=18$ &$g=19$ &$g=20$ \\\hline 
$q=0$ &11 &12 &13 &14 &15 &16 &17 &18 &19 &20 \\ 
$q=1$ &123 &174 &219 &291 &355 &453 &537 &666 &774 &936 \\ 
$q=2$ &176 &327 &538 &903 &1379 &2127 &3022 &4441 &5979 &8417 \\ 
$q=3$ &30 &75 &209 &448 &990 &1894 &3575 &6367 &10796 &17960 \\ 
$q=4$ &2 &3 &19 &34 &106 &295 &829 &1847 &4447 &9019 \\ 
$q=5$ &1 &1 &2 &2 &9 &18 &55 &116 &403 &986 \\ 
$q=6$ &&&1 &1 &2 &2 &7 &9 &36 &48 \\ 
$q=7$ &&&&&1 &1 &2 &2 &7 &7 \\ 
$q=8$ &&&&&&&1 &1 &2 &2 \\ 
$q=9$ &&&&&&&&&1 &1 \\ 
\hline
sum= &343 &592 &1001 &1693 &2857 &4806 &8045 &13467 &22464 &37396 \\\hline
\end{tabular}\\

\begin{tabular}{|c|c|c|c|c|c|c|c|c|c|c|}
\hline $\varrho_{g,q}$&$g=21$ &$g=22$ &$g=23$ &$g=24$ &$g=25$ &$g=26$ &$g=27$ &$g=28$ &$g=29$ &$g=30$ \\\hline 
$q=0$ &21 &22 &23 &24 &25 &26 &27 &28 &29 &30 \\ 
$q=1$ &1072 &1272 &1437 &1680 &1878 &2166 &2401 &2739 &3012 &3405 \\ 
$q=2$ &10966 &14826 &18774 &24770 &30539 &39321 &47697 &60083 &71711 &88938 \\ 
$q=3$ &28265 &44272 &66046 &99525 &140960 &204611 &281077 &394617 &525838 &720977 \\ 
$q=4$ &18673 &35178 &65533 &115252 &197836 &329568 &533479 &848091 &1304275 &2001344 \\ 
$q=5$ &2981 &7165 &17640 &37770 &84075 &166465 &331872 &615860 &1135074 &1989842 \\ 
$q=6$ &181 &464 &1383 &3603 &11141 &26864 &67991 &153882 &352322 &727680 \\ 
$q=7$ &25 &37 &94 &170 &652 &1679 &5300 &14899 &42738 &107050 \\ 
$q=8$ &7 &7 &23 &24 &85 &99 &321 &715 &2506 &7073 \\ 
$q=9$ &2 &2 &7 &7 &23 &23 &69 &83 &233 &331 \\ 
$q=10$ &1 &1 &2 &2 &7 &7 &23 &23 &68 &70 \\ 
$q=11$ &&&1 &1 &2 &2 &7 &7 &23 &23 \\ 
$q=12$ &&&&&1 &1 &2 &2 &7 &7 \\ 
$q=13$ &&&&&&&1 &1 &2 &2 \\ 
$q=14$ &&&&&&&&&1 &1 \\ 
\hline
sum= &62194 &103246 &170963 &282828 &467224 &770832 &1270267 &2091030 &3437839 &5646773 \\\hline
\end{tabular}\\

\end{tabular}}

    \end{center}}
  \caption{Number of semigroups of each genus and quasi-ordinarization number.}
  \label{t:ngq}
  \end{table}
}

\section{Analysis of $\varrho_{g,q}$}
\label{s:varrho}

Let us denote $o_{g,r}$ the number of numerical semigroups of genus $g$ and ordinarization number $r$ and $\varrho_{g,q}$ the number of numerical semigroups of genus $g$ and quasi-ordinarization number $r$

We can observe a behavior of $\varrho_{g,q}$ very similar to the behavior of 
$o_{g,r}$ introduced in \cite{Bras:ordinarization}.

Indeed, for odd $g$ and large $r$, it holds $\varrho_{g,q}=o_{g,r}$ 
and for even $g$ and large $q$, it holds $\varrho_{g,q}=o_{g,r+1}$. We will give a partial proof of these equalities at the end of this section.

It is conjcetured in \cite{Bras:ordinarization}
that, for each genus $g\in{\mathbb N_0}$ and each ordinarization number $r\in{\mathbb N_0}$, $$o_{g,r}\leq o_{{g+1},r}.$$
We can write the new conjecture below paralleling this.
\begin{Conjecture}
For each genus $g\in{\mathbb N_0}$ and each quasi-ordinarization number $q\in{\mathbb N_0}$, $$\varrho_{g,q}\leq \varrho_{{g+1},q}.$$
\end{Conjecture}

Now we will give some results on $\varrho_{g,q}$ for high quasi-ordinarization numbers.
First we will need Fre\u\i man's Theorem \cite{Freiman1,Freiman2} as formulated
in \cite{Nathanson}.
\begin{Theorem}[Fre{\u\i}man]
\label{t:freiman}
Let $A$ be a set of integers such that $\#A=k\geq 3$.
If $\#(A+A)\leq 3k-4$, then $A$ is a subset of an arithmetic 
progression of length $\#(A+A)-k+1\leq 2k-3$.
\end{Theorem}

The next lemma is a consequence of Fre\u\i man's Theorem. The lemma
tells that the first non-gaps of numerical semigroups of large quasi-ordinarization number 
must be even.

 \begin{Lemma}
 \label{lemma:parells}
 If a semigroup $\Lambda$ of genus $g$ has quasi-ordinarization number $q$ with
 $\frac{g+1}{3}\leq q\leq \lfloor\frac{g-1}{2}\rfloor$
 then all its non-gaps which are less than or equal to $g-1$ are even. 
 \end{Lemma}

\begin{proof}
  Suppose that $\Lambda$ is a semigroup of genus $g$ and quasi-ordinarization number $q\geq \frac{g+1}{3}$.

  Since the quasi-ordinarization is $q$, this means that $\lambda_0=0,\lambda_1,\dots,\lambda_q\leq g-1$ and $\lambda_{q+1}\geq g$, hence 
$\Lambda\cap[0,g-1]=\{\lambda_0,\lambda_1,\dots,\lambda_q\}$. 
Let $A=\Lambda\cap[0,g-1]$. By the previous equality, $\#A=q+1$. 
We have that the elements in $A+A$ are upper bounded by $2g-2$ and so 
$A+A\subseteq \Lambda\cap[0,2g-2]$.
Then $\#(A+A)\leq \#(\Lambda\cap[0,2g-2])< \#(\Lambda\cap[0,2g])$.
Since the Frobenius number of $\Lambda$ is at most $2g-1$,
$\#(\Lambda\cap[0,2g])=g+1$ and, so, $\#(A+A)\leq g$.
Now, since $q\geq\frac{g+1}{3}$ we have 
$g\leq 3q-1=3(q+1)-4=3(\#A)-4$ and we can apply Theorem~\ref{t:freiman} 
with $k=q+1$.
Thus we have that $A$ is a subset of an arithmetic progression of length at most $g-k+1=g-q$.

Let $d(A)$ be the difference between consecutive terms
of this arithmetic progression.
The number $d(A)$ can not be larger than or equal to three since otherwise
$\lambda_{q}\geq q\cdot d(A)\geq 3q\geq 3\frac{g+1}{3}>g$, a contradiction with $q$ being the quasi-ordinarization number.

If $d(A)=1$ then $A\subseteq[0,g-q-1]$ and
and so $\Lambda\cap[g-q,g-1]=\emptyset$.
We claim that 
in this case $A\subseteq\{0\}\cup[\lceil\frac{g}{2}\rceil,g-q-1]$.
Indeed, suppose that $x\in A$. Then $2x$ satisfies either $2x\leq g-q-1$ or $2x\geq g$.
If the second inequality is satisfied then it is obvious that 
$x\in\{0\}\cup[\lceil\frac{g}{2}\rceil,g-q-1]$.
If the first inequality is satisfied then we will prove that $mx\leq g-q-1$ for all $m\geq 2$ 
by induction on $m$ and this leads to $x=0$.
Indeed, if $mx\leq g-q-1$ then $x\leq \frac{g-q-1}{m}
\leq \frac{g-\frac{g+1}{3}-1}{m}
=\frac{2g-4}{3m}<\frac{2g}{3m}$.
Now $(m+1)x< \frac{2g(m+1)}{3m}=\frac{(2m+2)g}{3m}$ 
and since $m\geq 2$, we have $(m+1)x< \frac{(2m+m)g}{3m}=g$ and so
$(m+1)x\leq g-1$.
Since $(m+1)x$ is in $\Lambda\cap[0,g-1]=A\subseteq[0,g-q-1]$ this means that
$(m+1)x\leq g-q-1$ and this proves the claim.

Now, $A\subseteq\{0\}\cup[\lceil\frac{g}{2}\rceil,g-q-1]$
together with $\#A=q+1$
implies that 
$q\leq g-q-\lceil\frac{g}{2}\rceil=\lfloor\frac{g}{2}\rfloor-q\leq \frac{g}{2}-\frac{g+1}{3}=\frac{g-2}{6}<q$, a contradiction.

So, we deduce that $d(A)=2$, leading to the proof of the lemma.

\end{proof}

The next lemma was proved in \cite{Bras:ordinarization}.

\begin{Lemma}
\label{lemma:frob}
Suppose that a numerical 
semigroup $\Lambda$ has $\omega$ gaps between $1$ and $n-1$
and $n\geq 2\omega+2$ then 
\begin{enumerate}
\item $n\in\Lambda$,
\item the Frobenius number of $\Lambda$ is smaller than $n$,
\item the genus of $\Lambda$ is $\omega$.
\end{enumerate}
\end{Lemma}

Let $\Lambda$ be a numerical semigroup.
As in \cite{Bras:ordinarization}, let us say that a set
$B\subset{\mathbb N}_0$ is $\Lambda$-{\it closed}
if for any $b\in B$ and any $\lambda$ in $\Lambda$, the sum
$b+\lambda$ is either in $B$ or it is larger than $\max(B)$.
If $B$ is $\Lambda$-closed so is $B-\min(B)$. Indeed,
$b-\min(B)+\lambda=(b+\lambda)-\min(B)$ is either in $B-\min(B)$
or it is larger than $\max(B)-\min(B)=\max(B-\min(B))$.
The new $\Lambda$-closed set contains $0$.
We will denote by 
$C(\Lambda,i)$
the $\Lambda$-closed sets
of size $i$ that contain $0$. 

Let ${\mathcal S}_\gamma$ be the set of numerical semigroups of genus $\gamma$.
It was proved in \cite{Bras:ordinarization}
that, for $r$ an integer with $\frac{g+2}{3}\leq r \leq\lfloor\frac{g}{2}\rfloor$, it holds
$$o_{g,r}=\sum_{\Omega \in{\mathcal S}_{(\lfloor\frac{g}{2}\rfloor-r)}}\#C(\Omega,\left\lfloor\frac{g}{2}\right\rfloor-r+1).$$ 
We will see now that, for $q$ an integer with $\frac{g+1}{3}\leq q\leq \lfloor\frac{g-1}{2}\rfloor$, it holds  
$$\varrho_{g,q}=\sum_{\Omega\in{\mathcal S}_{(\lfloor\frac{g-1}{2}\rfloor-q)}}\#C(\Omega,\left\lfloor\frac{g-1}{2}\right\rfloor-q+1).$$ 
This proves that,
for $q$ an integer with $\frac{g+2}{3}\leq q\leq \lfloor\frac{g-1}{2}\rfloor$, we have $$\varrho_{g,q}=\left\{\begin{array}{ll}
o_{g,q}&\mbox{ if }g\mbox{ is odd,}\\
o_{g,q+1}&\mbox{ if }g\mbox{ is even.}\\
\end{array}\right.$$

\begin{Theorem}
\label{theorem:high}
 Let $g\in{\mathbb N}_0$, $g\geq 1$, and let $q$ be an integer with $\frac{g+1}{3}\leq q\leq \lfloor\frac{g-1}{2}\rfloor$. Define 
 $\omega=\lfloor\frac{g-1}{2}\rfloor-q$ 
 \begin{enumerate}
 \item
 If $\Omega$ is a numerical semigroup
 of genus $\omega$ and
 $B$ is a $\Omega$-closed set of size $\omega+1$ and first element equal to $0$ then  
 $$\{2j:j \in\Omega\}\cup \{2j-2\max(B)+2g+1:j\in B\}\cup(2g+{\mathbb N}_0)$$ 
 is a numerical semigroup of genus $g$ and quasi-ordinarization number equal to $q$.
 \item
 All numerical semigroups of genus $g$ and quasi-ordinarization number $q$ 
 can be uniquely written as
 $$\{2j:j \in\Omega\}\cup \{2j-2\max(B)+2g+1:j\in B\}\cup(2g+{\mathbb N}_0)$$ 
 for a unique numerical semigroup
 $\Omega$ of genus $\omega$ 
 and a unique $\Omega$-closed set $B$ of size $\omega+1$ and first element equal to $0$. 
 \item
 The number $\rho_{g,q}$ of numerical semigroups of genus $g$ and quasi-ordinarization number $q$
 depends only on $\omega$.
 It is exactly $$\sum_{\Omega\in{\mathcal S}_\omega}\#C(\Omega,\omega+1).$$
 \end{enumerate}
 \end{Theorem}

\begin{proof}
\begin{enumerate}
\item
Suppose that 
$\Omega$ is a numerical semigroup of genus $\omega$ 
and $B$ is a $\Omega$-closed set of size $\omega+1$ and first element equal to $0$. Let $X=\{2j:j \in\Omega\}$, $Y=\{2j-2\max(B)+2g+1:j\in B\}$, and $Z=(2g+{\mathbb N}_0)$.

First of all, let us see that the complement 
${\mathbb N}_0\setminus(X\cup Y\cup Z)$ 
has $g$ elements. Notice that all elements in $X$ are even while all elements in $Y$ are odd. So, $X$ and $Y$ do not intersect. 
Also the unique element in $Y\cap Z$ is $2g+1$. The number of elements in the complement will be 
\begin{eqnarray*}\#{\mathbb N}_0\setminus(X\cup Y\cup Z)&=&2g-\#\{x\in X:x<2g\}-\#Y+1\\&=&2g-\#\{s\in\Omega:s<g\}-\#B+1\\&=&2g-\omega-\#\{s\in\Omega:s<g\}.
\end{eqnarray*}

We know that all gaps of $\Omega$ are at most $2\omega-1=2(\lfloor\frac{g-1}{2}\rfloor-q)-1\leq g-1-2q-1<g$.
So, $\#\{s\in\Omega:s<g\}=g-\omega$ and we conclude that
$\#{\mathbb N}_0\setminus(X\cup Y\cup Z)=g$.

Before proving that $X\cup Y\cup Z$ 
is a numerical semigroup,
let us prove that the number of non-zero elements in $X\cup Y\cup Z$ 
which are smaller than or equal to $g-1$
is $q$. 
Once we prove that $X\cup Y\cup Z$ is a numerical semigroup,
this will mean, by Lemma~\ref{lemma:equivdef}, 
that it has quasi-ordinarization number $q$.
On the one hand, all elements in $Y$ are larger than $g-1$. 
Indeed, if $\lambda$ is the enumeration of $\Omega$ (i.e., $\Omega=\{\lambda_0,\lambda_1,\dots\}$ with $\lambda_i<\lambda_{i+1}$), then $\max(B)\leq\lambda_\omega\leq 2\omega=2\lfloor\frac{g-1}{2}\rfloor-2q\leq g-1-2\frac{g+1}{3}<\frac{g}{3} $. 
Now, for any $j\in B$, $2j-2\max(B)+2g+1>2g-2\max(B)>g$. 
On the other hand, 
all gaps of $\Omega$ are at most 
$2\omega-1=2\lfloor\frac{g-1}{2}\rfloor-2q-1< g-\frac{2(g+1)}{3}-1<\frac{g}{3}-1$ and so
all the even integers not belonging to $X$ are less than $g$.
So, the number of non-zero non-gaps of $X\cup Y\cup Z$ smaller than or equal to
$g-1$ is $\lfloor\frac{g-1}{2}\rfloor-\omega=q$.

To see that $X\cup Y\cup Z$ 
is a numerical semigroup we only need to see that
it is closed under addition. It is obvious that $X+Z\subseteq Z$, $Y+Z\subseteq Z$, $Z+Z\subseteq Z$. It is also obvious that $X+X\subseteq X$ since $\Omega$ is a numerical semigroup and that $Y+Y\subseteq Z$ since, as we proved before, all elements in $Y$ are larger than $g$.

It remains to see that 
$X+Y\subseteq X\cup Y\cup Z$.
Suppose that $x\in X$ and $y\in Y$. Then $x=2i$ for some $i\in\Omega$ and $y=2j-2\max(B)+2g+1$ 
for some $j\in B$. Then $x+y=2(i+j)-2\max(B)+2g+1$. Since $B$ is $\Omega$-closed, we have that either 
$i+j\in B$ and so $x+y\in Y$ or $i+j>\max(B)$. In this case $x+y\in Z$. So, $X+Y\subseteq Y\cup Z$.

\item
First of all notice that,
since the Frobenius number of a semigroup $\Lambda$ of genus $g$ is
smaller than $2g$, it holds
$$\Lambda\cap(2g+{\mathbb N}_0)=(2g+{\mathbb N}_0).$$
For any numerical semigroup $\Lambda$,
the set $\Omega=\{\frac{j}{2}:j\in\Lambda\cap(2{\mathbb N}_0)\}$ is a numerical semigroup.
If $\Lambda$ has quai-ordinarization number $q\geq\frac{g+1}{3}$ then, 
by Lemma~\ref{lemma:parells},
$$\Lambda\cap[0,g-1]=(2\Omega)\cap[0,g-1].$$
The semigroup $\Omega$ has exactly
$q+1$ non-gaps between $0$ and $\lfloor\frac{g-1}{2}\rfloor$ and
$\omega=\lfloor\frac{g-1}{2}\rfloor-q$
gaps between $0$ and $\lfloor\frac{g-1}{2}\rfloor$.
We can use Lemma~\ref{lemma:frob} with
$n=\lfloor\frac{g+1}{2}\rfloor$
since
$$2\omega+2=2\left\lfloor\frac{g-1}{2}\right\rfloor-2q+2\leq g-1-\frac{2(g+1)}{3}+2=\frac{g+1}{3},$$
which implies $2\omega+2\leq \frac{g+1}{3}\leq \lfloor\frac{g+1}{2}\rfloor=n$.
Then the genus of $\Omega$ is $\omega$ and the Frobenius number of $\Omega$ is 
at most 
$\lfloor\frac{g+1}{2}\rfloor$.
This means in particular that all even integers larger than $g-1$
belong to $\Lambda$.

Define $D=(\Lambda\cap[0,2g])\setminus2\Omega$.
That is, $D$ is the set of odd non-gaps of $\Lambda$ smaller than $2g$.
We claim that 
$\bar B=\{\frac{j-1}{2}: j\in D\cup\{2g+1\}\}$
is a $\Omega$-closed set of size $\omega+1$.
The size follows from the fact that the number of non-gaps of $\Lambda$ between 
$g$ and $2g$ is $g-q$ and that the number of even integers 
in the same interval is
$\lceil\frac{g+1}{2}\rceil$.
Suppose that $\lambda\in \Omega$ and $b\in \bar B$.
Then $b=\frac{j-1}{2}$ for some $j$ in $D\cup\{2g+1\}$ and $b+\lambda=
\frac{(j+2\lambda)-1}{2}$.
If $\frac{(j+2\lambda)-1}{2}\geq\max(\bar B)=\frac{(2g+1)-1}{2}$ we are done. Otherwise we have $j+2\lambda\leq 2g$. Since $\Lambda$ is a numerical semigroup and
both $j,2\lambda\in\Lambda$, it holds $j+2\lambda\in\Lambda\cap[0,2g]$.
Furthermore, $j+2\lambda$ is odd since so is $j$. So, $b+\lambda$ is either
larger than $\max(\bar B)$ or it is in $\bar B$.
Then $B=\bar B-\min(\bar B)$ is a $\Lambda$-closed set of size $\omega+1$ and 
first element zero.

\item
    The previous two points define a bijection between
    the set of numerical semigroups in ${\mathcal S}_g$ of quasi-ordinarization number $q$  and the set $\{C(\Omega,\omega+1): \Omega\in{\mathcal S}_\omega\}.$
    Hence,
    $\rho_{g,q}=\sum_{\Omega\in{\mathcal S}_\omega}\#C(\Omega,\omega+1).$

\end{enumerate}
\end{proof}

\begin{Corollary}
Suppose that $\frac{g+2}{3}\leq q\leq \lfloor\frac{g-1}{2}\rfloor$. 
Then,
$$\varrho_{g,q}=\left\{\begin{array}{ll}
o_{g,q}&\mbox{ if }g\mbox{ is odd,}\\
o_{g,q+1}&\mbox{ if }g\mbox{ is even.}\\
\end{array}\right.$$
\end{Corollary}

 Define, as in \cite{Bras:ordinarization}, the sequence $f_\omega$ by 
 $f_\omega=\sum_{\Omega\in{\mathcal S}_\omega}\#C(\Omega,\omega+1).$
 The first elements in the sequence, from $f_0$ to $f_{\maxomega}$ are
 \begin{center}
 \resizebox{.9\textwidth}{!}{$\begin{array}{|c|cccccccccccccccc|}
 \hline
 \omega & 0 & 1 & 2 & 3 & 4 & 5 & 6 &  7 & 8 & 9 & 10 & 11 & 12 & 13 & 14 & 15  \\
 \hline
 f_\omega & 1 & 2 & 7& 23& 68& 200& 615& 1764& 5060& 14626& 41785& 117573& 332475& 933891& 2609832& 7278512 \\\hline\end{array}$}
 \end{center}

 We remark that this sequence appears in \cite{bernardini}, where Bernardini and Torres proved that the number of numerical semigroups of genus $3\omega$ whose number of even gaps is $\omega$ is exactly $f_\omega$. It corresponds to the entry \href{http://oeis.org/A210581}{A210581} of The On-Line Encyclopedia of Integer Sequences \cite{oeis}.

 We can deduce the values $\varrho_{g,q}$ using the values in the previous table together with Theorem~\ref{theorem:high} 
 for any $g$, whenever $q\geq\max(\frac{g+2}{3},\lfloor\frac{g-1}{2}\rfloor-\maxomega)$.

 The next corollary is a consequence of the fact that the sequence $f_\omega$ is increasing for $\omega$ between $0$ and $\maxomega$.
 \begin{Corollary}
 For any $g\in{\mathbb N}$ and 
 any $q\geq\max(\frac{g}{3}+1,\lfloor\frac{g}{2}\rfloor-\maxomega)$, it holds $\varrho_{g,q}\geq \varrho_{g+1,q}$.
 \end{Corollary}

 If we proved that $f_\omega\leq f_{\omega+1}$ for any $\omega$, 
 this would imply $\varrho_{g,q}\leq \varrho_{g+1,q}$ for any $q>\frac{g}{3}$.

 \section{The forest ${\mathscr F}_g$}
 \label{s:forest}
Fix a genus $g$. We can define a graph in which the nodes are all semigroups of that genus and whose edges connect each semigroup to its quasi-ordinarization transform, if it is not itself.
The graph is a forest ${\mathscr F}_g$
rooted at all ordinary and quasi-ordinary semigroups of genus $g$.
In particular, the quasi-ordinarization transform defines, 
for each fixed genus and conductor, a tree rooted at the unique 
quasi-ordinary semigroup of that genus and conductor, given in Lemma~\ref{l:roots}.
See ${\mathscr F}_4$ in Figure~\ref{fquatre},
${\mathscr F}_6$ in Figure~\ref{fsis}, and
${\mathscr F}_7$ in Figure~\ref{fset}.

\begin{figure}
  \begin{adjustbox}{max totalsize={.35\textwidth}{.35\textheight},left}\begin{tikzpicture}[grow=right, every node/.style = {align=left}]\tikzset{level 1+/.style={level distance=12cm}}\Tree[.\nongap{0}\gap{1}\gap{2}\gap{3}\gap{4}\nongap{5}\nongap{6}\nongap{7}\nongap{8} ]
 \end{tikzpicture}\end{adjustbox}

\bigskip

\begin{adjustbox}{max totalsize={.9\textwidth}{.9\textheight},left}\begin{tikzpicture}[grow=right, every node/.style = {align=left}]\tikzset{level 1+/.style={level distance=12cm}}\Tree[.{$\nongap{0}\gap{1}\gap{2}\gap{3}\nongap{4}\gap{5}\nongap{6}\nongap{7}\nongap{8}\dots $} [.{$\nongap{0}\gap{1}\gap{2}\nongap{3}\gap{4}\gap{5}\nongap{6}\nongap{7}\nongap{8}\dots $}  ]
 ]
 \end{tikzpicture}\end{adjustbox}

\bigskip

\begin{adjustbox}{max totalsize={.9\textwidth}{.9\textheight},left}\begin{tikzpicture}[grow=right, every node/.style = {align=left}]\tikzset{level 1+/.style={level distance=12cm}}\Tree[.{$\nongap{0}\gap{1}\gap{2}\gap{3}\nongap{4}\nongap{5}\gap{6}\nongap{7}\nongap{8}\dots $}  ]
 \end{tikzpicture}\end{adjustbox}

\bigskip

\begin{adjustbox}{max totalsize={.9\textwidth}{.9\textheight},left}\begin{tikzpicture}[grow=right, every node/.style = {align=left}]\tikzset{level 1+/.style={level distance=12cm}}\Tree[.{$\nongap{0}\gap{1}\gap{2}\gap{3}\nongap{4}\nongap{5}\nongap{6}\gap{7}\nongap{8}\dots $} [.{$\nongap{0}\gap{1}\gap{2}\nongap{3}\gap{4}\nongap{5}\nongap{6}\gap{7}\nongap{8}\dots $}  ]
[.{$\nongap{0}\gap{1}\nongap{2}\gap{3}\nongap{4}\gap{5}\nongap{6}\gap{7}\nongap{8}\dots $}  ]
 ]
 \end{tikzpicture}\end{adjustbox}

  \caption{${\mathscr F}_4$}
\label{fquatre}
\end{figure}

\begin{figure}
  \begin{adjustbox}{max totalsize={.5\textwidth}{.5\textheight},left}\begin{tikzpicture}[grow=right, every node/.style = {align=left}]\tikzset{level 1+/.style={level distance=12cm}}\Tree[.\nongap{0}\gap{1}\gap{2}\gap{3}\gap{4}\gap{5}\gap{6}\nongap{7}\nongap{8}\nongap{9}\nongap{10}\nongap{11}\nongap{12} ]
 \end{tikzpicture}\end{adjustbox}

\bigskip

\begin{adjustbox}{max totalsize={\textwidth}{\textheight},left}\begin{tikzpicture}[grow=right, every node/.style = {align=left}]\tikzset{level 1+/.style={level distance=12cm}}\Tree[.{$\nongap{0}\gap{1}\gap{2}\gap{3}\gap{4}\gap{5}\nongap{6}\gap{7}\nongap{8}\nongap{9}\nongap{10}\nongap{11}\nongap{12}\dots $} [.{$\nongap{0}\gap{1}\gap{2}\gap{3}\nongap{4}\gap{5}\gap{6}\gap{7}\nongap{8}\nongap{9}\nongap{10}\nongap{11}\nongap{12}\dots $}  ]
[.{$\nongap{0}\gap{1}\gap{2}\gap{3}\gap{4}\nongap{5}\gap{6}\gap{7}\nongap{8}\nongap{9}\nongap{10}\nongap{11}\nongap{12}\dots $}  ]
 ]
 \end{tikzpicture}\end{adjustbox}

\bigskip

\begin{adjustbox}{max totalsize={\textwidth}{\textheight},left}\begin{tikzpicture}[grow=right, every node/.style = {align=left}]\tikzset{level 1+/.style={level distance=12cm}}\Tree[.{$\nongap{0}\gap{1}\gap{2}\gap{3}\gap{4}\gap{5}\nongap{6}\nongap{7}\gap{8}\nongap{9}\nongap{10}\nongap{11}\nongap{12}\dots $} [.{$\nongap{0}\gap{1}\gap{2}\gap{3}\gap{4}\nongap{5}\gap{6}\nongap{7}\gap{8}\nongap{9}\nongap{10}\nongap{11}\nongap{12}\dots $}  ]
[.{$\nongap{0}\gap{1}\gap{2}\nongap{3}\gap{4}\gap{5}\nongap{6}\gap{7}\gap{8}\nongap{9}\nongap{10}\nongap{11}\nongap{12}\dots $}  ]
[.{$\nongap{0}\gap{1}\gap{2}\gap{3}\gap{4}\nongap{5}\nongap{6}\gap{7}\gap{8}\nongap{9}\nongap{10}\nongap{11}\nongap{12}\dots $}  ]
 ]
 \end{tikzpicture}\end{adjustbox}

\bigskip

\begin{adjustbox}{max totalsize={\textwidth}{\textheight},left}\begin{tikzpicture}[grow=right, every node/.style = {align=left}]\tikzset{level 1+/.style={level distance=12cm}}\Tree[.{$\nongap{0}\gap{1}\gap{2}\gap{3}\gap{4}\gap{5}\nongap{6}\nongap{7}\nongap{8}\gap{9}\nongap{10}\nongap{11}\nongap{12}\dots $} [.{$\nongap{0}\gap{1}\gap{2}\gap{3}\nongap{4}\gap{5}\gap{6}\nongap{7}\nongap{8}\gap{9}\nongap{10}\nongap{11}\nongap{12}\dots $}  ]
[.{$\nongap{0}\gap{1}\gap{2}\gap{3}\gap{4}\nongap{5}\gap{6}\nongap{7}\nongap{8}\gap{9}\nongap{10}\nongap{11}\nongap{12}\dots $}  ]
[.{$\nongap{0}\gap{1}\gap{2}\gap{3}\nongap{4}\gap{5}\nongap{6}\gap{7}\nongap{8}\gap{9}\nongap{10}\nongap{11}\nongap{12}\dots $}  ]
[.{$\nongap{0}\gap{1}\gap{2}\gap{3}\gap{4}\nongap{5}\nongap{6}\gap{7}\nongap{8}\gap{9}\nongap{10}\nongap{11}\nongap{12}\dots $}  ]
[.{$\nongap{0}\gap{1}\gap{2}\gap{3}\gap{4}\nongap{5}\nongap{6}\nongap{7}\gap{8}\gap{9}\nongap{10}\nongap{11}\nongap{12}\dots $}  ]
 ]
 \end{tikzpicture}\end{adjustbox}

\bigskip

\begin{adjustbox}{max totalsize={\textwidth}{\textheight},left}\begin{tikzpicture}[grow=right, every node/.style = {align=left}]\tikzset{level 1+/.style={level distance=12cm}}\Tree[.{$\nongap{0}\gap{1}\gap{2}\gap{3}\gap{4}\gap{5}\nongap{6}\nongap{7}\nongap{8}\nongap{9}\gap{10}\nongap{11}\nongap{12}\dots $} [.{$\nongap{0}\gap{1}\gap{2}\gap{3}\nongap{4}\gap{5}\gap{6}\nongap{7}\nongap{8}\nongap{9}\gap{10}\nongap{11}\nongap{12}\dots $}  ]
[.{$\nongap{0}\gap{1}\gap{2}\nongap{3}\gap{4}\gap{5}\nongap{6}\gap{7}\nongap{8}\nongap{9}\gap{10}\nongap{11}\nongap{12}\dots $}  ]
 ]
 \end{tikzpicture}\end{adjustbox}

\bigskip

\begin{adjustbox}{max totalsize={\textwidth}{\textheight},left}\begin{tikzpicture}[grow=right, every node/.style = {align=left}]\tikzset{level 1+/.style={level distance=12cm}}\Tree[.{$\nongap{0}\gap{1}\gap{2}\gap{3}\gap{4}\gap{5}\nongap{6}\nongap{7}\nongap{8}\nongap{9}\nongap{10}\gap{11}\nongap{12}\dots $} [.{$\nongap{0}\gap{1}\gap{2}\gap{3}\gap{4}\nongap{5}\gap{6}\nongap{7}\nongap{8}\nongap{9}\nongap{10}\gap{11}\nongap{12}\dots $} [.{$\nongap{0}\gap{1}\gap{2}\gap{3}\nongap{4}\nongap{5}\gap{6}\gap{7}\nongap{8}\nongap{9}\nongap{10}\gap{11}\nongap{12}\dots $}  ]
 ]
[.{$\nongap{0}\gap{1}\gap{2}\gap{3}\nongap{4}\gap{5}\nongap{6}\gap{7}\nongap{8}\nongap{9}\nongap{10}\gap{11}\nongap{12}\dots $} [.{$\nongap{0}\gap{1}\nongap{2}\gap{3}\nongap{4}\gap{5}\nongap{6}\gap{7}\nongap{8}\gap{9}\nongap{10}\gap{11}\nongap{12}\dots $}  ]
 ]
[.{$\nongap{0}\gap{1}\gap{2}\nongap{3}\gap{4}\gap{5}\nongap{6}\nongap{7}\gap{8}\nongap{9}\nongap{10}\gap{11}\nongap{12}\dots $}  ]
 ]
 \end{tikzpicture}\end{adjustbox}

  \caption{${\mathscr F}_6$}
\label{fsis}
\end{figure}

\begin{figure}
  \begin{adjustbox}{max totalsize={.35\textwidth}{.35\textheight},left}\begin{tikzpicture}[grow=right, every node/.style = {align=left}]\tikzset{level 1+/.style={level distance=12cm}}\Tree[.\nongap{0}\gap{1}\gap{2}\gap{3}\gap{4}\gap{5}\gap{6}\gap{7}\nongap{8}\nongap{9}\nongap{10}\nongap{11}\nongap{12}\nongap{13}\nongap{14} ]
 \end{tikzpicture}\end{adjustbox}

\bigskip

\begin{adjustbox}{max totalsize={.9\textwidth}{.9\textheight},left}\begin{tikzpicture}[grow=right, every node/.style = {align=left}]\tikzset{level 1+/.style={level distance=12cm}}\Tree[.{$\nongap{0}\gap{1}\gap{2}\gap{3}\gap{4}\gap{5}\gap{6}\nongap{7}\gap{8}\nongap{9}\nongap{10}\nongap{11}\nongap{12}\nongap{13}\nongap{14}\dots $} [.{$\nongap{0}\gap{1}\gap{2}\gap{3}\gap{4}\nongap{5}\gap{6}\gap{7}\gap{8}\nongap{9}\nongap{10}\nongap{11}\nongap{12}\nongap{13}\nongap{14}\dots $}  ]
[.{$\nongap{0}\gap{1}\gap{2}\gap{3}\gap{4}\gap{5}\nongap{6}\gap{7}\gap{8}\nongap{9}\nongap{10}\nongap{11}\nongap{12}\nongap{13}\nongap{14}\dots $}  ]
 ]
 \end{tikzpicture}\end{adjustbox}

\bigskip

\begin{adjustbox}{max totalsize={.9\textwidth}{.9\textheight},left}\begin{tikzpicture}[grow=right, every node/.style = {align=left}]\tikzset{level 1+/.style={level distance=12cm}}\Tree[.{$\nongap{0}\gap{1}\gap{2}\gap{3}\gap{4}\gap{5}\gap{6}\nongap{7}\nongap{8}\gap{9}\nongap{10}\nongap{11}\nongap{12}\nongap{13}\nongap{14}\dots $} [.{$\nongap{0}\gap{1}\gap{2}\gap{3}\nongap{4}\gap{5}\gap{6}\gap{7}\nongap{8}\gap{9}\nongap{10}\nongap{11}\nongap{12}\nongap{13}\nongap{14}\dots $}  ]
[.{$\nongap{0}\gap{1}\gap{2}\gap{3}\gap{4}\nongap{5}\gap{6}\gap{7}\nongap{8}\gap{9}\nongap{10}\nongap{11}\nongap{12}\nongap{13}\nongap{14}\dots $}  ]
[.{$\nongap{0}\gap{1}\gap{2}\gap{3}\gap{4}\gap{5}\nongap{6}\gap{7}\nongap{8}\gap{9}\nongap{10}\nongap{11}\nongap{12}\nongap{13}\nongap{14}\dots $} [.{$\nongap{0}\gap{1}\gap{2}\gap{3}\gap{4}\nongap{5}\nongap{6}\gap{7}\gap{8}\gap{9}\nongap{10}\nongap{11}\nongap{12}\nongap{13}\nongap{14}\dots $}  ]
 ]
[.{$\nongap{0}\gap{1}\gap{2}\gap{3}\gap{4}\nongap{5}\gap{6}\nongap{7}\gap{8}\gap{9}\nongap{10}\nongap{11}\nongap{12}\nongap{13}\nongap{14}\dots $}  ]
[.{$\nongap{0}\gap{1}\gap{2}\gap{3}\gap{4}\gap{5}\nongap{6}\nongap{7}\gap{8}\gap{9}\nongap{10}\nongap{11}\nongap{12}\nongap{13}\nongap{14}\dots $}  ]
 ]
 \end{tikzpicture}\end{adjustbox}

\bigskip

\begin{adjustbox}{max totalsize={.9\textwidth}{.9\textheight},left}\begin{tikzpicture}[grow=right, every node/.style = {align=left}]\tikzset{level 1+/.style={level distance=12cm}}\Tree[.{$\nongap{0}\gap{1}\gap{2}\gap{3}\gap{4}\gap{5}\gap{6}\nongap{7}\nongap{8}\nongap{9}\gap{10}\nongap{11}\nongap{12}\nongap{13}\nongap{14}\dots $} [.{$\nongap{0}\gap{1}\gap{2}\gap{3}\nongap{4}\gap{5}\gap{6}\gap{7}\nongap{8}\nongap{9}\gap{10}\nongap{11}\nongap{12}\nongap{13}\nongap{14}\dots $}  ]
[.{$\nongap{0}\gap{1}\gap{2}\gap{3}\gap{4}\gap{5}\nongap{6}\gap{7}\nongap{8}\nongap{9}\gap{10}\nongap{11}\nongap{12}\nongap{13}\nongap{14}\dots $} [.{$\nongap{0}\gap{1}\gap{2}\nongap{3}\gap{4}\gap{5}\nongap{6}\gap{7}\gap{8}\nongap{9}\gap{10}\nongap{11}\nongap{12}\nongap{13}\nongap{14}\dots $}  ]
 ]
[.{$\nongap{0}\gap{1}\gap{2}\gap{3}\gap{4}\gap{5}\nongap{6}\nongap{7}\gap{8}\nongap{9}\gap{10}\nongap{11}\nongap{12}\nongap{13}\nongap{14}\dots $}  ]
[.{$\nongap{0}\gap{1}\gap{2}\gap{3}\nongap{4}\gap{5}\gap{6}\nongap{7}\nongap{8}\gap{9}\gap{10}\nongap{11}\nongap{12}\nongap{13}\nongap{14}\dots $}  ]
[.{$\nongap{0}\gap{1}\gap{2}\gap{3}\gap{4}\gap{5}\nongap{6}\nongap{7}\nongap{8}\gap{9}\gap{10}\nongap{11}\nongap{12}\nongap{13}\nongap{14}\dots $}  ]
 ]
 \end{tikzpicture}\end{adjustbox}

\bigskip

\begin{adjustbox}{max totalsize={.9\textwidth}{.9\textheight},left}\begin{tikzpicture}[grow=right, every node/.style = {align=left}]\tikzset{level 1+/.style={level distance=12cm}}\Tree[.{$\nongap{0}\gap{1}\gap{2}\gap{3}\gap{4}\gap{5}\gap{6}\nongap{7}\nongap{8}\nongap{9}\nongap{10}\gap{11}\nongap{12}\nongap{13}\nongap{14}\dots $} [.{$\nongap{0}\gap{1}\gap{2}\gap{3}\nongap{4}\gap{5}\gap{6}\gap{7}\nongap{8}\nongap{9}\nongap{10}\gap{11}\nongap{12}\nongap{13}\nongap{14}\dots $}  ]
[.{$\nongap{0}\gap{1}\gap{2}\gap{3}\gap{4}\nongap{5}\gap{6}\gap{7}\nongap{8}\nongap{9}\nongap{10}\gap{11}\nongap{12}\nongap{13}\nongap{14}\dots $}  ]
[.{$\nongap{0}\gap{1}\gap{2}\gap{3}\gap{4}\gap{5}\nongap{6}\gap{7}\nongap{8}\nongap{9}\nongap{10}\gap{11}\nongap{12}\nongap{13}\nongap{14}\dots $} [.{$\nongap{0}\gap{1}\gap{2}\nongap{3}\gap{4}\gap{5}\nongap{6}\gap{7}\gap{8}\nongap{9}\nongap{10}\gap{11}\nongap{12}\nongap{13}\nongap{14}\dots $}  ]
[.{$\nongap{0}\gap{1}\gap{2}\gap{3}\nongap{4}\gap{5}\nongap{6}\gap{7}\nongap{8}\gap{9}\nongap{10}\gap{11}\nongap{12}\nongap{13}\nongap{14}\dots $}  ]
 ]
[.{$\nongap{0}\gap{1}\gap{2}\gap{3}\gap{4}\nongap{5}\gap{6}\nongap{7}\gap{8}\nongap{9}\nongap{10}\gap{11}\nongap{12}\nongap{13}\nongap{14}\dots $}  ]
[.{$\nongap{0}\gap{1}\gap{2}\gap{3}\gap{4}\gap{5}\nongap{6}\nongap{7}\gap{8}\nongap{9}\nongap{10}\gap{11}\nongap{12}\nongap{13}\nongap{14}\dots $}  ]
[.{$\nongap{0}\gap{1}\gap{2}\gap{3}\gap{4}\nongap{5}\gap{6}\nongap{7}\nongap{8}\gap{9}\nongap{10}\gap{11}\nongap{12}\nongap{13}\nongap{14}\dots $}  ]
[.{$\nongap{0}\gap{1}\gap{2}\gap{3}\gap{4}\gap{5}\nongap{6}\nongap{7}\nongap{8}\gap{9}\nongap{10}\gap{11}\nongap{12}\nongap{13}\nongap{14}\dots $}  ]
[.{$\nongap{0}\gap{1}\gap{2}\gap{3}\gap{4}\gap{5}\nongap{6}\nongap{7}\nongap{8}\nongap{9}\gap{10}\gap{11}\nongap{12}\nongap{13}\nongap{14}\dots $}  ]
 ]
 \end{tikzpicture}\end{adjustbox}

\bigskip

\begin{adjustbox}{max totalsize={.9\textwidth}{.9\textheight},left}\begin{tikzpicture}[grow=right, every node/.style = {align=left}]\tikzset{level 1+/.style={level distance=12cm}}\Tree[.{$\nongap{0}\gap{1}\gap{2}\gap{3}\gap{4}\gap{5}\gap{6}\nongap{7}\nongap{8}\nongap{9}\nongap{10}\nongap{11}\gap{12}\nongap{13}\nongap{14}\dots $} [.{$\nongap{0}\gap{1}\gap{2}\gap{3}\gap{4}\nongap{5}\gap{6}\gap{7}\nongap{8}\nongap{9}\nongap{10}\nongap{11}\gap{12}\nongap{13}\nongap{14}\dots $}  ]
 ]
 \end{tikzpicture}\end{adjustbox}

\bigskip

\begin{adjustbox}{max totalsize={.9\textwidth}{.9\textheight},left}\begin{tikzpicture}[grow=right, every node/.style = {align=left}]\tikzset{level 1+/.style={level distance=12cm}}\Tree[.{$\nongap{0}\gap{1}\gap{2}\gap{3}\gap{4}\gap{5}\gap{6}\nongap{7}\nongap{8}\nongap{9}\nongap{10}\nongap{11}\nongap{12}\gap{13}\nongap{14}\dots $} [.{$\nongap{0}\gap{1}\gap{2}\gap{3}\gap{4}\gap{5}\nongap{6}\gap{7}\nongap{8}\nongap{9}\nongap{10}\nongap{11}\nongap{12}\gap{13}\nongap{14}\dots $} [.{$\nongap{0}\gap{1}\gap{2}\gap{3}\gap{4}\nongap{5}\nongap{6}\gap{7}\gap{8}\nongap{9}\nongap{10}\nongap{11}\nongap{12}\gap{13}\nongap{14}\dots $}  ]
[.{$\nongap{0}\gap{1}\gap{2}\gap{3}\nongap{4}\gap{5}\nongap{6}\gap{7}\nongap{8}\gap{9}\nongap{10}\nongap{11}\nongap{12}\gap{13}\nongap{14}\dots $} [.{$\nongap{0}\gap{1}\nongap{2}\gap{3}\nongap{4}\gap{5}\nongap{6}\gap{7}\nongap{8}\gap{9}\nongap{10}\gap{11}\nongap{12}\gap{13}\nongap{14}\dots $}  ]
 ]
[.{$\nongap{0}\gap{1}\gap{2}\nongap{3}\gap{4}\gap{5}\nongap{6}\gap{7}\nongap{8}\nongap{9}\gap{10}\nongap{11}\nongap{12}\gap{13}\nongap{14}\dots $}  ]
 ]
[.{$\nongap{0}\gap{1}\gap{2}\gap{3}\gap{4}\nongap{5}\gap{6}\nongap{7}\gap{8}\nongap{9}\nongap{10}\nongap{11}\nongap{12}\gap{13}\nongap{14}\dots $}  ]
[.{$\nongap{0}\gap{1}\gap{2}\gap{3}\nongap{4}\gap{5}\gap{6}\nongap{7}\nongap{8}\gap{9}\nongap{10}\nongap{11}\nongap{12}\gap{13}\nongap{14}\dots $}  ]
 ]
 \end{tikzpicture}\end{adjustbox}

  \caption{${\mathscr F}_7$}
\label{fset}
\end{figure}

In the forest ${\mathscr F}_g$ we know that the parent of a numerical semigroup that is not a root is its quasi-ordinarization transform.
Let us analyze now, what are the children of a numerical semigroup. The next result is well know and can be found, for instance, in \cite{NS}. We use $\Lambda^*$ to denote $\Lambda\setminus\{0\}$.

\begin{Lemma}\label{l:fundamentalgaps}
Suppose that $\Lambda$ is a numerical semigroup and that $a\in{\mathbb N}_0\setminus \Lambda$. The set $a\cup \Lambda$ is a numerical semigroup if and only if
\begin{itemize}
\item $a+\Lambda^*\subseteq \Lambda^*$,
\item $2a\in \Lambda$,
\item $3a\in \Lambda$.
\end{itemize}
\end{Lemma}

The elements $a\in{\mathbb N}_0\setminus \Lambda$ such that 
$a+\Lambda\subseteq \Lambda$, are denoted {\it pseudo-Frobenius numbers} of $\Lambda$.
The elements $a\in{\mathbb N}_0\setminus \Lambda$ such that 
$\{2a,3a\}\subseteq \Lambda$, are denoted {\it fundamental gaps} of $\Lambda$.
The elements satisfying the three conditions will be called {\it candidates}.

Suppose that a numerical semigroup $\Lambda$ with Frobenius number $F$ has children in ${\mathscr F}_g$.
Let $e_1,\dots,e_r$ be the generators of $\Lambda$ between the subconductor and $F-1$.
For $i=1,\dots,r$, let $c^i_1,\dots,c^i_{k_{i}}$ be the candidates of $\Lambda\setminus\{e_i\}$. The children of $\Lambda$ in ${\mathscr F}_g$
are the semigroups of the form $\Lambda\setminus\{e_i\}\cup\{c^i_j\}$, for $i=1,\dots,r$ and $j=1,\dots,k_i$.

\section{Relating ${\mathscr F}_g$, ${\mathscr T}_g$, and ${\mathscr T}$}
\label{s:Ts}

Now we analyze the relation between the kinship of different nodes in ${\mathscr F}_g$, ${\mathscr T}_g$, and ${\mathscr T}$.
If two semigroups are children of the same semigroup $\Lambda$, then they are called \emph{siblings}. If $\Lambda_1$ and $\Lambda_2$ are siblings, and $\Lambda_3$ is a child of $\Lambda_2$, then we say that $\Lambda_3$ is a \emph{niece$/$nephew} of $\Lambda_1$. 

Let $q(\Lambda)$ denote the quasi-ordinarization of $\Lambda$.
The next lemmas are quite immediate from the definitions.

\begin{Lemma}
If $\Lambda_1$ is a child of $\Lambda_2$ in $\mathscr T$ then 
$q(\Lambda_1)$ is a niece/nephew of $q(\Lambda_2)$ in $\mathscr T$. 
\end{Lemma}

As an example, $\Lambda_1=\{0,4,5,8,9,10,12,\dots\}$
is a child of $\Lambda_2=\{0,4,5,8,\dots\}$ in ${\mathscr T}$, while $q(\Lambda_1)=\{0,5,7,8,9,10,12,\dots\}$ is a niece of
$q(\Lambda_2)=\{0,5,6,8,\dots\}$ in $\mathscr T$. 

\begin{Lemma}
If $\Lambda_1$ and $\Lambda_2$ are siblings 
in $\mathscr T$ then they are siblings
in ${\mathscr T}_g$ but not in
${\mathscr F}_g$.
\end{Lemma}

As an example, $\Lambda_1=\{0,5,7,9,10,11,12,14,\dots\}$ and $\Lambda_2=\{0,5,7,9,10,12,\dots\}$ are siblings in ${\mathscr T}$ and in ${\mathscr T}_7$ (see Figure~\ref{genere}),
but they are not siblings in
${\mathscr F}_7$ (see Figure~\ref{fset}).

\begin{Lemma}
If $\Lambda_1$ and $\Lambda_2$ are siblings
in ${\mathscr T}_g$ 
then $q(\Lambda_1)$ and $q(\Lambda_2)$ are siblings in 
$\mathscr T$.
\end{Lemma}

As an example, $\Lambda_1=\{0,3,6,9,10,12,\dots\}$ and $\Lambda_2=\{0,5,6,10,\dots\}$ are siblings in ${\mathscr T}_7$ (see Figure~\ref{genere}), and
$q(\Lambda_1)=\{0,6,8,9,10,12,\dots\}$
 and $q(\Lambda_2)=\{0,6,8,10,\dots\}$ are siblings in 
$\mathscr T$.
 
As a consequence of the previous two lemmas, we get this last lemma.

\begin{Lemma}
If $\Lambda_1$ and $\Lambda_2$ are siblings 
in $\mathscr T$ then
$q(\Lambda_1)$ and $q(\Lambda_2)$ are siblings in
$\mathscr T$.
\end{Lemma}

As an example, $\Lambda_1=\{0,5,7,9,10,11,12,14,\dots\}$ and $\Lambda_2=\{0,5,7,9,10,12,\dots\}$ are siblings in ${\mathscr T}$ and
$q(\Lambda_1)=\{0,7,8,9,10,11,12,14,\dots\}$ and $q(\Lambda_2)=\{0,7,8,9,10,12,\dots\}$ are siblings in ${\mathscr T}$.

\section{Conclusion}

Quasi-ordinary semigroups are those semigroups that have all gaps except one in a row, while ordinary semigroups have all gaps in a row.

We defined a quasi-ordinarization transform that, applied repeatedly to a non-ordinary numerical semigroup stabilizes in a quasi-ordinary semigroup of the same genus.

From this transform, fixing a genus $g$, we can define a forest ${\mathscr F}_g$ whose nodes are all semigroups of genus $g$, whose roots are all ordinary and quasi-ordinary semigroups of that genus, and whose edges connect each non-ordinary and non-quasi-ordinary numerical semigroup to its quasi-ordinarization transform.

We conjectured that the number of numerical semigroups in 
${\mathscr F}_g$ at a given depth 
is at most the number of numerical semigroups in 
${\mathscr F}_{g+1}$ at the same depth.
We provided a proof of the conjecture for the largest possible depths. Proving this conjecture for all depths, would prove the conjecture that $n_{g+1}\geq n_{g}$. Hence, we expect our work to be a step forward towards the proof of the conjectured increasingness of the sequence $n_g$.


\footnotesize 

\end{document}